\documentclass[a4paper, 12pt]{article}
\usepackage{theorem}
\usepackage[dvips]{graphics}
\newtheorem{definition}{Definition}
\newtheorem{proposition}{Lemma}
\newtheorem{theorem}{Theorem}
\newtheorem{cor}{Corollary}
\usepackage{graphics, amssymb, amsmath}
\def\disp{\displaystyle}

\def\C{\mathcal{C}}
\def\mb{\mathbb}
\def\mc{\mathcal}

\def\lan{\langle}
\def\ran{\rangle}
\begin{document}
\centerline{\large{An extended existence result for quadratic BSDEs with jumps}}
 \centerline{\large{with application to the utility maximization problem}}
%opening\title{Efficient hedging and risk minimization}
\vspace{0.5cm}
\centerline{Marie-Amelie Morlais\footnote{A large part of the content of this work is in my PhDthesis defended at the university of Rennes 1 in October 2007 and supervised by Professor Ying Hu}}
\vspace{0.2cm}
\centerline{Universit\'e du Maine}
\vspace{0.2cm}
\centerline{Avenue Olivier Messaien, 72085 Le Mans}

%%\footnote{This author is sponsored by Credit Suisse. The major part of this work was done at the ETH of Zurich, Switzerland}
%%\centerline{Ramistrasse 101, 8006 Zurich,} Tel: +41 44 632 5859,
\centerline{E-mail: marieamelie.morlais@free.fr
}
\vspace{0.5cm}

%Comparison result and uniqueness as a by product
%+ A priori estimates (fundamental)\\

\begin{abstract}
 In this study, we consider the exponential utility maximization problem in the context of a jump-diffusion model. To solve this problem, we rely on the dynamic programming principle and we derive from it a quadratic BSDE with jumps. Since this quadratic BSDE\footnote{The notation of quadratic BSDE refers to the growth with respect of the variable $z$ of the generator $f:(s,z,u ) \to f(s,z, u)$.} is driven both by a Wiener process and a Poisson random measure having a Levy measure with infinite mass, our main work consists in establishing a new existence result for the specific BSDE introduced.
\end{abstract}

\section{Introduction}
In this paper, our motivation is to study the exponential utility maximization problem with portfolio constraints in the context of a discontinuous filtration. To handle this optimization problem, which is formulated at any time under a conditional form, the approach consists in using both the martingale optimality principle and BSDE techniques: this approach is the same as in the previous papers \cite{Becherer}, \cite{ManiaetSchw} and \cite{morlais2} already dealing with the same problem. However and contrary to the papers \cite{Becherer} or \cite{morlais2} already dealing with a discontinuous model, the originality of the present paper is that we study existence for a specific class of quadratic BSDEs with jumps without assuming the finiteness of the Levy measure. Relaxing this last hypothesis, we have to establish a new existence result for the BSDE already introduced in \cite{morlais2}, which is the main achievement of this paper. Concerning the financial problem under study, the main objectives are the characterization of the value process
in terms of the solution of an explicit BSDE as well as the characterization of optimal strategies. 
\\
\indent  To obtain the main result, that is the existence of solutions of the specific BSDE introduced by using the dynamic programming principle, we first define an auxiliary BSDE (more precisely, we introduce a new generator which is explicitely given in terms of the first one) and we then prove the existence result for the auxiliary BSDE under an additional constraint on the norm of the bounded terminal condition. For the general case, i.e. when considering a BSDE whose terminal condition is an arbitrary bounded random variable, we provide an explicit construction. In a last step, we first establish a correspondence result between solutions of the auxiliary BSDE and those of the original one and we then prove existence of a solution of the original BSDE for any arbitrary random variable. In a last section, we come back and solve the original financial problem. 
%Finally and in a last step, we come back to the original financial problem.
\\
\indent The present paper is structured as follows: in Section 2, we describe the financial model and we give preliminary notations. Then, in Sections 3 and 4, we state and prove the main results for the BSDE introduced in Section 2.
Last section consists in using results of the two previous sections to provide answers to the original financial problem. Lengthy proofs are relegated to the appendix.
\section{The model and preliminaries}
\hspace{0.5cm} We consider a probability space ($\Omega$, $\mathbb{F}$, $\mathbb{P}$) equipped with two independent stochastic processes:
\begin{enumerate}
\item[.] A standard (one dimensional) brownian motion: $W$ =$(W_{t})_{t \in [0,T]}$.
\item[.] A real-valued
Poisson point process $p$ defined on $[0, T] \times \mb{R}\setminus{\{0\}}$. Referring to chapter 2 in \cite{Ikedawatanabe}, we denote by $N_{p}(ds,dx)$ the associated counting measure, whose compensator is assumed to be of the form  $$\disp{\hat{N}_{p}(ds,dx) = n(dx)ds}.$$
 $n(dx)$ (also denoted by $n$ in the sequel) stands for the Levy measure which is positive and satisfies
$$n(\{0\})= 0  \;\textrm{and}  \; \disp{\int_{\mb{R}\setminus{\{0\}}}(1 \wedge |x|)^{2}n(dx) } < \infty.$$
%\begin{equation}\label{eq: conditionmeslevy}
%\big( n(\{0\}) =0 \big) \;\textrm{and} \; \left( n((1 \wedge |x|)^{2}) := \disp{\int_{\mb{R}\setminus{\{0\}}}(1 %\wedge |x|)^{2}n(dx) < \infty}\right). \end{equation}
 \end{enumerate}
 These two processes $W$ and $\tilde{N}_{p}$ are considered on $[0, T]$, where $T$ stands for the horizon or maturity time in the financial context and, in all the sequel, $T$ is assumed to be fixed and deterministic.
%In that paper, we relax the assumption of finiteness of the measure $n$, which is the one used both in %\cite{Becherer} or in \cite{morlais2}. 
We also denote by
 $\mc{F}$ the filtration generated by the two processes $W$ and $N_{p}$ (and completed by $\mathcal{N}$, consisting in all the $\mb{P}$-null sets). 
Using the same notations as in \cite{Ikedawatanabe}, we denote by $\tilde{N}_{p}(ds,dx)$ ($\tilde{N}_{p}(ds, dx) := N_{p}(ds, dx) - \hat{N}_{p}(ds, dx)$) the compensated measure, which is a martingale random measure: in particular, for any predictable and locally square integrable process $K$,
the stochastic integral $K \cdot \tilde{N}_{p} := \disp{\int K_{s}(x) \tilde{N}_{p}(ds,dx)}$ is a locally square integrable martingale. \\
We denote by $Z \cdot W $ (resp. $U \cdot \tilde{N}_{p}$) the stochastic integral of $Z$ w.r.t. $W$ (resp. the stochastic integral of $U$ w.r.t. $\tilde{N}_{p}$). Since the filtration $\mc{F}$ has the predictable representation property, then, for any local martingale $M$ of $\mc{F}$, there exists two predictable processes $Z$ and $U$ such that  
$$\forall \; t, \quad M_{t} = M_{0} + \big(Z \cdot W \big)_{t} + \big( U \cdot \tilde{N}_{p}\big)_{t}. $$ 
(In Section 2.2, we provide a definition of the Hilbert spaces, where these stochastic integrals are considered). In all the paper, we will make use of the notation $|\cdot|_{\infty}$ to refer to the norm in $L^{\infty}(\mc{F}_{T})$ of any bounded $\mc{F}_{T}$-measurable random variable.

 \subsection{Preliminaries about BSDEs}
In the sequel, we denote by $\mc{S}^{\infty}(\mb{R})$ the set of all adapted processes $Y$ with c\`adl\`ag paths (c\`adl\`ag stands for right continuous with left limits) such that $$\textrm{ess} \disp{\sup_{t, \omega}}|Y_{t}(\omega)| < \infty,$$
and, for any $p$, $p > 0$, we denote by $\mc{S}^{p}$ the set of c\`adl\`ag processes $Y$ such that
$$\disp{\mb{E}\left(\sup_{t} |Y_{t}|^{p} \right)} < \infty. $$
 We also introduce the set $L^{2}(W)$ consisting of all predictable processes $Z$ such that 
%the stochastic integral $Z \cdot W $ satisfies
$$\disp{\mathbb{E}\left(\int_{0}^{T}|Z_{s}|^{2}ds \right) <\infty. } $$
  and the set $L^{2}(\tilde{N}_{p})$ consisting of all $\mathcal{P} \otimes \mathcal{B}(\mathbb{R}\setminus\{0\}) $-measurable processes $U$ such that 
$$ \disp{ \mathbb{E}\left(\int_{[0, T]\times \mathbb{R}\setminus\{0\}}|U_{s}(x)|^{2}n(dx)ds \right) < \infty.}$$
$\mathcal{P}$ stands for the $\sigma$-field of all predictable sets of $[0, T] \times \Omega$ and $\mathcal{B}(\mathbb{R}\setminus\{0\})$ the Borel field of $\mb{R}\setminus\{0\}$.
The set $L^{0}(n)$, which is also denoted by $L^{0}(n, \mb{R}, \mb{R}\setminus\{0\})$ in \cite{Becherer}, consists of all the functions $u$ mapping $\mb{R}$ in $\mb{R}\setminus\{0\}$ and it is equipped with the topology of convergence in measure. Finally, $L^{2}(n)$ stands for the subset of all functions in $L^{0}(n)$ such that: $\disp{\mb{E}\big(\int_{0}^{T}|u(x)|^{2}n(dx)\big)} < \infty$ and $L^{\infty}(n)$ stands for the subset of all functions $u $ in $L^{0}(n) $ which takes bounded values (almost surely).
  \\
\indent A solution of a BSDE with jumps of the form 
\begin{equation}\label{eq: EDSRavecsauts} 
Y_{t} = B + \int_{t}^{T} f(s, Y_{s-}, Z_{s}, U_{s})ds - \int_{t}^{T} Z_{s}dW_{s} - \int_{t}^{T}\int_{\mathbb{R}^{*}}U_{s}(x)\tilde{N}_{p}(ds, dx),
\end{equation}
which is characterized by a bounded terminal condition $B$ and a generator $f$ satisfying 
$$\disp{\int_{0}^{T}|f(s, Y_{s}, Z_{s}, U_{s})|ds} < \infty, \; \mb{P}\textrm{-a.s}.,$$
is a triple of processes ($Y$, $Z$, $U$) which is in
$ S^{\infty}(\mb{R}) \times L^{2}(W) \times L^{2}(\tilde{N}_{p} )$.  
In this paper, we study a specific class of BSDE with jumps of the previous form.
Besides and since we do not work on a brownian filtration, the processes $Z$ and $U$ have to be predictable, for any solution of the BSDE (\ref{eq: EDSRavecsauts}) .\\
%\indent A solution of such BSDEs with jumps is usually defined on
%$ S^{2} \times L^{2}(W) \times L^{2}(\tilde{N}_{p} )$, $S^{2}$ being equipped with the following norm: $|Y|_{S^{2}} = %\mathbb{E}(\disp{\sup_{t \in [0, T]}|Y_{t}|^{2}})^{\frac{1}{2}}$ (our main references for studies on BSDEs with %jumps, are \cite{BarlesBuck} and \cite{Royer2}). But, in this paper, the results of the aforementionned papers cannot %be applied, since the generator of the BSDE we are interested in, does not satisfy the usual conditions (it is not %Lipschitz w.r.t. its variable $z$). 

 \subsection{Description of the model}
For sake of completeness, we provide the description of the financial context which is similar as in \cite{morlais2}.
The financial market consists in one risk-free asset (assumed to have zero interest rate) and one single risky asset, whose price process is denoted by $S$.
More precisely, the stock price process is a one dimensional semimartingale satisfying 
\begin{eqnarray}\label{eq: eqavecsauts} dS_{s} = S_{s-}\left( b_{s}ds + \sigma_{s}dW_{s} + \int_{\mathbb{R}^{*}}{\beta_{s}(x)\tilde{N}_{p}(ds,dx)} \right) .  \end{eqnarray}
All processes $b$, $\sigma$ and $\beta$ are assumed to be bounded and predictable and, in addition, $\beta$ satisfies: $\beta > -1$. This last condition implies that the stochastic exponential $\mc{E}(\beta \cdot \tilde{N}_{p})$ is positive, $\mb{P}$-a.s.: hence, the price process $S$ is itself almost surely positive.
The boundedness of $\beta$, $\sigma$ and $\theta$ ensures both existence and uniqueness results for the SDE (\ref{eq: eqavecsauts}).
Then, provided that: $\sigma \neq 0$, we can define $\theta$ by: $\theta_{s} =\sigma_{s}^{-1}b_{s}$ ($\mathbb{P}$-a.s. and for all $s$).  
 The process $\theta$, also called market price of risk process, is supposed to be bounded and, under this assumption, the measure $\mb{P}^{\theta}$ with density
 $$ \dfrac{d\mathbb{P}^{\theta}}{d\mathbb{P}} = \mathcal{E}_{T}(-\int_{0}^{.}\theta_{s}dW_{s}), $$
 is a risk-neutral measure, which means that, under $\mb{P}^{\theta}$, the price process $S$ is a local martingale. \\
\indent In what follows, we introduce the usual notions of trading strategies and self financing portfolio, assuming that all trading strategies are constrained to take their values in a closed set denoted by $\mc{C}$. In a first step and to make easier the proofs, this set $\mc{C}$ is supposed to be compact\footnote{As in \cite{morlais2}, the compactness assumption on the constraint $\mc{C}$ ensures that the BMO properties given in ($H_{2}$) in Section 3.1 are satisfied: thanks to these properties, we can prove a comparison result for the BSDE with generator having the generator defined in (\ref{eq: generateurf}). In a last section of this aforementionned paper and by means of an approximating procedure, the existence result is obtained without this restrictive hypothesis.}. Due to the presence of constraints in this model with finite horizon $T$, not any $\mc{F}_{T}$-measurable random variable $B$ is attainable by using contrained strategies. In that context, we adress the problem of characterizing dynamically the value process associated to the exponential utility maximization problem (in the sequel, we denote by $U_{\alpha}$ the exponential utility function with parameter $\alpha$, which is defined on $\mb{R}$ by: $ U_{\alpha}(\cdot) = - \exp(- \alpha \cdot)$). 
\begin{definition}\label{tradingstrat}
A predictable $\mathbb{R}$-valued process $\pi$ is a self-financing trading strategy, if it takes its values in a constraint set $\mc{C}$ and if the process $X^{\pi, t, x}$ such that    \begin{equation}\label{eq: wealthproc}
\forall \; s \in [t,\;T], \quad X_{s}^{\pi, t, x}:= \disp{x + \int_{t}^{s}\pi_{s}\frac{dS_{s}}{S_{s-}}},
\end{equation}
 is in the space $\mc{H}^{2}$ of semimartingales (see chapter 4, \cite{Protter}).
Such a process $X^{\pi}=X^{\pi, t, x}$ stands for the wealth of an agent having strategy $\pi$ and wealth $x$ at time $t$.
 \end{definition}
\indent Now, as soon as the constraint set $\mc{C}$ is compact, the set consisting of all constrained strategies satisfies an additional integrability property.
\begin{proposition}\label{classequality}
 Under the assumption of compactness of the constraint set $\mathcal{C}$,
all trading strategies $\pi:=(\pi_{s})_{s \in [t, T]}$ as introduced in Definition \ref{tradingstrat} satisfy 
 \begin{equation}\label{eq: uniformint}   
\{\disp{ \exp(-\alpha X_{\tau}^{\pi}), \tau \; \mathcal{F}\textrm{-stopping time} \;  } \} \; \textrm{is a uniformly integrable family}.
 \end{equation}
\end{proposition}
For the proof of this lemma, we refer to \cite{morlais2}.
%\begin{proposition}\label{classequality}
% Under the assumption of compactness of the constraint set $\mathcal{C}$,
%all trading strategies $\pi:=(\pi_{s})_{s \in [t, T]}$ as introduced in Definition \ref{tradingstrat} satisfy 
% \begin{equation}\label{eq: integrability}   
%\disp{\exists \; \epsilon > 0, \;  \quad\mb{E}\big(\exp(|(\alpha + \epsilon) X_{t}^{\pi}| )\big) < \infty},
 %\end{equation}
%or equivalently, this means that the process ($-\alpha X^{\pi}$) admits an exponential moment of order greater than 1.
%\end{proposition}
We make use of the notation $\mathcal{A}_{t}$ for the admissibility set (in the case when $t=0$, we simply denote it by $\mathcal{A}$.): in this notation, the subscript $t$ indicates that we start the wealth dynamics at time $t$: more precisely, this set consists in all the strategies whose restriction to the interval $ [0, t] $ is equal to zero and which satisfy both Definition \ref{tradingstrat} and the condition (\ref{eq: uniformint}). This last integrability condition is of great use in Section 4 to justify the expression of the value process (and, more particularly, to justify the supermartingale property of some family of processes as already introduced in \cite{ImkelleretHu} in a Brownian setting). 
% which has already been used in \cite{ImkelleretHu} in a Brownian setting. The usual and much more restrictive admissibility condition consists in assuming that the wealth process $X^{\pi}$ is bounded from below (uniformly over all strategies $\pi$).\\
 \indent To conclude this paragraph, we introduce the notion of BMO martingales which can also be found in \cite{Dellacheriemeyer}: a martingale $M$ is said to be in the class of BMO martingales if there exists a constant $c$, $c > 0$, such that, for all $\mc{F}$-stopping time $\tau$,
$$   \textrm{ess}\displaystyle{\sup_{\Omega}\mathbb{E}^{\mathcal{F}_{\tau}}(\lan M \ran _{T}-\lan M \ran_{\tau} )} \leq c^{2}
\;\textrm{and} \; |\Delta M_{\tau}|^{2} \leq c^{2}.   $$
(In the continuous case, the BMO property follows from the first condition, whereas, in the discontinuous setting, we need to ensure the boundedness of the jumps of $M$).
 The following result, referred as Kazamaki's criterion and also stated in \cite{Kazamaki}, relates the martingale property of a stochastic exponential to a BMO property.

\begin{proposition}\label{kamazaki} \textbf{(Kazamaki's criterion)} $\quad$Let $\delta$ be such that: 0 $< \delta < \infty $ and $M$ a BMO martingale satisfying: $\Delta M_{t} \geq -1 + \delta$, $\mathbb{P}$-a.s. and for all $t$, then $\mc{E}(M)$ is a true martingale.
\end{proposition}

\section{The quadratic BSDE with jumps}
\subsection{Main assumptions}
In all the sequel, we use the explicit form of the generator $f$
\begin{eqnarray}\label{eq: generateurf}
 \disp{ f(s,z,u) = \displaystyle{\inf_{\pi \in
 \mathcal{C}}\left(\frac{\alpha}{2}|\pi\sigma_{s} - (z+ \frac{\theta_{s}}{\alpha})|^{2} + |u - \pi \beta_{s}|_ {\alpha}\right) -
  \theta_{s} z - \frac{|\theta_{s}|^{2}}{2\alpha}}},\\
\nonumber \end{eqnarray}
where the processes $\beta$, $\theta$ and $\sigma$ are defined in Section 2.1.
This expression of the generator will be justified in Section 4. 
We introduce the notation $|\cdot|_{\alpha}$ as being the convex functional such that 
\[ \begin{array}{ll}
\forall \; u \in (L^{2} \cap L^{\infty})(n),\; |u|_{\alpha} & =\disp{\int_{\mb{R}\setminus\{0\}}\frac{\exp(\alpha u(x)) - \alpha u(x) -1}{\alpha}n(dx)}, \\
\\
& =\disp{\int_{\mb{R}\setminus\{0\}}g_{\alpha}(u(x)) n(dx) },\\ 
\end{array}\] 
with the real function $g_{\alpha}$ defined by: $g_{\alpha}(y) = \frac{\exp(\alpha y) - \alpha y - 1 }{\alpha}$.
In all the paper, $B$ is a bounded $\mc{F}_{T}$-measurable random variable and we use these two standing assumptions on the generator $f$
\begin{enumerate}
 \item[$(H_{1})$.] 
The first assumption denoted by ($H_{1}$) consists in specifying both a lower and an upper bound for $f$
\[  \begin{array}{l}
 \forall \; z, u \in \mb{R} \times (L^{2} \cap L^{\infty})(n) \\
\\
- \theta_{s} z - \frac{|\theta_{s}|^{2}}{2\alpha}  \le f(s, z, u) \le \frac{\alpha}{2}|z|^{2}  + |u|_{\alpha} ,\; \mb{P}\textrm{-a.s. and for all}\;s. \\
\end{array} \]
\item[($H_{2}$).] The second assumption, referred as ($H_{2}$), consists in two estimates: 
the first one deals with the increments of the generator $f$ w.r.t. $z$
\[  \begin{array}{l}
  \exists \; C > 0, \;\kappa \in BMO(W), \;\forall \; z, \; z' \in \mathbb{R},\; \forall u \in L^{2}(n(dx)),   \\
\\
|f(s, z, u) - f(s, z', u )| \leq  C(\kappa_{s} + |z| + |z'|)|z - z'| \\
\\
\end{array} \]
The second estimate deals with the increments w.r.t. $u$
\[  \begin{array}{l}
      \forall  z \in \mathbb{R}, \;\forall \;  u, u' \in (L^{2}\cap L^{\infty})(n(dx)),     \\
\\
f(s, z, u) - f(s, z, u') \leq \disp{\int_{\mb{R}\setminus{\{0\}}}\gamma_{s}(u, u')(u(x) - u'(x))n(dx)},\;  \\
    \end{array} \]
with the following expression for $\gamma_{s}(u, u^{'})$ for all $s$
\[ \begin{array}{l}
\disp{ \; \gamma_{s}(u, u')  \; =  }\\
\\
\displaystyle{\sup_{\pi \in \mathcal{C}} \left(\int_{0}^{1} g_{\alpha}^{'}(\lambda(u- \pi \beta_{s}) + (1 - \lambda)(u'- \pi\beta_{s})(x))d\lambda \right)}\mathbf{1}_{u \geq u'} \quad \quad \quad  \\
  \\
 \; \; + \;  \displaystyle{ \inf_{\pi \in \mathcal{C}} \left(\int_{0}^{1}g_{\alpha}^{'}(\lambda(u- \pi \beta_{s}) + (1 - \lambda)(u'- \pi \beta_{s})(x)d\lambda \right)}\mathbf{1}_{u < u'},\quad \quad \quad  \\
 \end{array}  \]
and this last expression holds, for any fixed
 $s, \omega$.
Considering now two arbitrary predictable processes $U$, $U^{'}$ taking their values in $ L^{2} \cap L^{\infty} (n)$ and if we define the process
 $\tilde{\gamma}$
 for all $s$ by 
\begin{equation}\label{eq: processgamma}
\tilde{\gamma}_{s} =
\gamma_{s}(U_{s}, U^{'}_{s}), \end{equation}
then, $\tilde{\gamma}$ is a predictable process and it is explicitely given in terms of both the predictable processes $U$, $ U^{'}$ and $\beta$.
 For the proof of these two estimates and the justification of the expression of $\gamma$, the reader is referred to \cite{morlais2}. 
To conclude this paragraph, 
we justify the BMO property of the process given by (\ref{eq: processgamma}): for this, we
use the compactness of $\mc{C}$ and we assume that both processes $U$ and $U'$ take their values $ L^{2} \cap L^{\infty}(n(dx))$ and that: $ |U_{s}|_{L^{\infty}(n)}, |U_{s}'|_{L^{\infty}(n)} \leq K$, to argue that
$$ \exists \; \delta_{K}, \bar{C}_{K} > 0,\; \textrm{s.t.} \quad  - 1 + \delta_{K} \leq \gamma_{s}(U_{s}, U_{s}^{'}) \leq \bar{C}_{K}, $$
which entails, in particular, that this process is in BMO($\tilde{N}_{p}$).
We rely on this BMO property in the proof of the uniqueness result to justify the use of Girsanov's theorem.

\end{enumerate}
 \subsection{Theoretical results} 
 To prove the main existence result, which is the existence of solutions of BSDEs with generator $f$ given by (\ref{eq: generateurf}) and terminal condition $B$ ($B$ being an arbitrary bounded random variable), we need to consider an auxiliary BSDE with parameters ($\tilde{f}, \tilde{B}$): more precisely, we consider the generator $\tilde{f}$ defined in terms of $f$ as follows
$$ \tilde{f}(s, z, u) = f(s, z -\frac{\theta_{s}}{\alpha}, u) - f(s,-\frac{\theta_{s}}{\alpha}, 0).$$
 In the first step, we motivate the introduction of this auxiliary BSDE
 by proving an existence result: to do this, the idea consists in establishing precise a priori estimates given by (\ref{eq: estimessent}) to justify, in a second step, a new stability result, which is similar as in \cite{morlais2}. This will be done under an explicit constraint on the terminal condition.
In the following theorem, we state the two main existence results of this paper.
%We state here the main result of that paper which gives existence for a solution of the BSDE($f, \; \bar{B}$) with %$f$ given by (\ref{eq: generateurf}) and the terminal condition $\bar{B}$ in $\mc{S}^{p}$ for any $p$, $p > 0$ (the %expression of $\bar{B}$ will be explicitely characterized in terms of a bounded random variable $B$ and of the %parameters $\theta$ and $\alpha$ introduced in Section 2). 
\begin{theorem}\label{existence}
(i) For any BSDE of the form (\ref{eq: EDSRavecsauts}) with generator $\tilde{f}$ and terminal condition $B$ satisfying 
$$\forall \; \; k >0, \quad \quad \mb{E}\left(\exp(k |B|) \right) < \infty,
 $$
%and such that $B$ is a bounded random variable,
 there exists at least one solution ($Y, Z, U$) such that $\exp(Y)$ is in $\mc{S}^{p}$, for any $p$, $p > 0$, and ($Z, U$) is in $ L^{2}(W) \times L^{2}(\tilde{N}_{p})$.\\
(ii) For any BSDE of the form (\ref{eq: EDSRavecsauts}) with generator $f$ and terminal condition $\bar{B}$, such that $\bar{B}$ is an arbitrary bounded random variable, there exists at least one solution ($\bar{Y}, \; \bar{Z}, \;\bar{U}$) in $ \mc{S}^{\infty} \times L^{2}(W) \times L^{2}(\tilde{N}_{p})$. 
\end{theorem}
%To obtain this result, we first justify the existence result for an auxiliary quadratic BSDE with jumps (with bounded %terminal condition). 
For later use, we provide here some a priori estimates for solutions of BSDEs with jumps having a bounded terminal condition (the proof of this lemma can be found in \cite{morlais2}).

	\begin{proposition}\label{estim2}
For any BSDE of the form (\ref{eq: EDSRavecsauts}) with a generator $g$ satisfying ($H_{1}$) and a bounded terminal condition $B$, there exists three explicit constants $C_{1}$, $C_{2}$ and $C_{3}$ given in terms of $|B|_{\infty}$, $|\theta|_{S^{\infty}(\mb{R})}$ and $\alpha$, and such that, for any solution ($Y$, $Z$, $U$) in $S^{\infty}(\mb{R}) \times L^{2}(W) \times L^{2}(\tilde{N}_{p})$ and for any $\mathcal{F}$-stopping time $\tau$, $\tau$ taking its values in $[0, T]$,
	\[       \begin{array}{l}  
	    (i) \; \mb{P}\textrm{-a.s. and for all} \;\; t, \; t \in [0, T], \;\; C_{1} \leq
	     Y_{t} 
	     \leq C_{2},\;   \\
\\
	     (ii) \;   \mathbb{E}^{\mathcal{F}_{\tau}}\disp(\int_{\tau}^{T}|Z_{s}|^{2}ds + \int_{\tau}^{T}{\int_{\mathbb{R}^{*}}{|U_{s}(x)|^{2}
n(dx)}ds} \disp) \leq C_{3}.   \\

	     \nonumber \end{array}  \]
	       \end{proposition}

\begin{cor}\label{equivalence}
Under the same assumptions than in Lemma \ref{estim2} on the parameters $g$ and $B$ and for any solution ($Y$, $Z$, $U$) in $S^{\infty}(\mb{R}) \times L^{2}(W) \times L^{2}(\tilde{N}_{p})$ of the BSDE (\ref{eq: EDSRavecsauts}), 
\begin{itemize}
 \item there exists a predictable version $\tilde{U}$ of $U$ such that: $\tilde{U} \equiv U $ (in $L^{2}(\tilde{N}_{p})$). Noting $U$ instead of $\tilde{U}$, this process satisfies\footnote{Here and contrary to Corollary 1 in \cite{morlais2}, since the Levy measure satisfies: $n(\mb{R}^{*}) = \infty$, we cannot deduce that $u$ takes its values in $L^{2}(n)$, using the fact that it is in $L^{\infty}(n)$.}
$$  |U_{s}|_{L^{\infty}(n)} \leq 2 |Y|_{S^{\infty}(\mb{R})}. $$
\item The following equivalence result
\begin{eqnarray}\label{eq: relationeq}
\nonumber \exists \; C \; > 0, \;  \; \frac{1}{C} \mathbb{E}\int_{[0, T] \times \mb{R}\setminus{\{0\}}}|U_{s}(x)|^{2}n(dx)ds  \leq \disp{\mathbb{E}\int_{0}^{T}|U_{s}|_{\alpha}ds} \\
\quad \quad \quad \quad \leq \disp{C \mathbb{E}\int_{[0, T] \times \mb{R}\setminus{\{0\}}}|U_{s}(x)|^{2}n(dx)ds,} \end{eqnarray}
holds for a constant $C$ depending only on $\alpha$ and $|Y|_{S^{\infty}(\mb{R})}$.
\end{itemize}
\end{cor} 
 \subsection{ Proof of the main existence result}
First and for sake of clarity, we give an outline of the content of this section. To prove Theorem \ref{existence}, we proceed with the following steps \\
\textbullet $\;$  In a first step, we introduce the auxiliary generator $\tilde{f}$ such that
\begin{equation}\label{eq: generateurauxiliaire}
 \tilde{f}(s, z, u) = f(s, z -\frac{\theta_{s}}{\alpha}, u) - f(s,-\frac{\theta_{s}}{\alpha}, 0), 
\end{equation}
and we then establish an existence result for the BSDEs given by ($\tilde{f}, \;\frac{B}{N} $) by providing a sufficient condition on the integer $N$.\\
\textbullet $\;$  In a second step and to prove existence for the BSDE given by ($\tilde{f}, \; B$) for any bounded $\mc{F}_{T}$-measurable random variable $B$, we proceed with an iterative \footnote{The construction is iterative in the following sense that the generator $f^{i+1}$ is defined in terms of $f^{i}$.} procedure. To this end, we construct a sequence of BSDEs given by ($f^{i}, \;\frac{B}{N}$) such that,
% in particular, it satisfies: $\disp{\sum_{i=1}^{k} f^{i}} = \tilde{f}$, for any $k$, $k \ge 1$.
 under the assumption that there exists a solution ($\tilde{Y}^{i}, \tilde{Z}^{i}, \tilde{U}^{i}$) up to step $k$, the triple $(\bar{Y}^{k}, \bar{Z}^{k}, \bar{U}^{k})$ with: $\bar{Y}^{k} = \disp{\sum_{i}\tilde{Y}^{i}}$, solves the BSDE with parameters ($\tilde{f}, \;\displaystyle{\sum_{i=1}^{k}\frac{B}{N}}$).
Provided this construction can be iterated up to step $N$, the process $Y$ defined by: $Y = \bar{Y}^{N}$ solves the BSDE with parameters ($\tilde{f}, \; B$).\\
\textbullet $\;$ The third step consists in establishing a correspondence result between a solution of the BSDE given by the parameters ($\tilde{f} , \; B$) and a solution of the BSDE with parameters ($f, \bar{B}$), with $\bar{B}$ explicitely given in terms of $B$.\\
\textbullet $\;$ Finally, in a last step, we extend the results of Step 2 to the case when the terminal condition may be unbounded (but admits at least exponential moments of any order). This is done by using the same methodology as in \cite{BriandetHu}: this step allows to prove existence for solutions of the BSDE with generator $f$ when the terminal condition is arbitrary and bounded.\\ 
\subsubsection{Step 1: first approximation}
\paragraph*{Construction and basic properties}
Since we are dealing with a BSDE with jumps whose generator has quadratic growth, we rely on the same procedure as in \cite{morlais2}: this consists in constructing an approximating sequence of generators denoted by $(f^{m})$. To this end, we introduce the constant  
 $M$, the truncation function $\rho_{m}$ and the measure $n^{m}$ as follows\\
(i) $M = 2(C_{1} + C_{2})$ (these two constants are given in (i)(a), Lemma \ref{estim2}).\\
 (ii) $\rho_{m}$ is an arbitrary truncation function at least continuously differentiable and such that: $\rho_{m}(z) = 0$, if $|z| \ge m+1$ and
 $\rho_{m}(z) = 1$, if $|z| \le m$, and $0 \le \rho_{m}(z) \le 1$ if $ 0 \le  z \le 1$.\\
(iii) $n^{m}$ is the finite measure defined by
$$ n^{m}(dx) = \mathbf{1}_{|x| \ge \frac{1}{m}}n(dx).
$$
This being set, we define the sequence
 ($f^{m}$) by
\[  \begin{array}{ll}
f^{m}(s, z, u) = & \disp{\inf_{ \pi \in \mc{C}} \left(\frac{\alpha}{2}|\pi\sigma_{s} - (z+ \frac{\theta_{s}}{\alpha}) |^{2}\rho_{m}(z) + \int_{\mb{R}^{*}}g_{\alpha}(u - \pi\beta_{s})\rho_{M}(u(x)) n^{m}(dx)\right)}\\
  & \quad  - z\theta_{s} - \frac{|\theta_{s}|^{2}}{2\alpha},   \end{array} \]
 and we then introduce $(f^{1,m})$ by setting
$$f^{1, m}(s, z , u) = f^{m}(s, z- \frac{\theta_{s}}{\alpha}, u) - f(s,- \frac{\theta_{s}}{\alpha}, 0).  $$
Since $0$ is in the set $\mc{C}$, the infimum in the expression of $f^{m}(s, \frac{- \theta_{s}}{\alpha}, 0)$ is equal to zero and hence, we obtain: 
$ f^{m}(s, \frac{- \theta_{s}}{\alpha}, 0) = f(s, \frac{- \theta_{s}}{\alpha}, 0) = \frac{|\theta_{s}|^{2}}{\alpha},$ implying that
$$ \forall \; m,  \quad f^{1, m}(s, 0 , 0)  \equiv 0, \quad \mb{P}\textrm{-a.s. and for all} \; s. $$
We provide below a list of the essential properties satisfied by ($f^{1, m}$)
\begin{enumerate}
\item[1.] Due to the truncation procedure, the generator $f^{1, m}$ is lipschitz with respect to $z$ and $u$, i.e.
there exists a constant $C_{m}$ depending only on the bounded parameters $\theta$, $\beta$, and on the constants
$\alpha$
and $\disp{\sup_{\pi \in
\mc{C}|}|\pi|}$,
such that
$$ |f^{1, m}(s, z, u) - f^{1, m}(s, z^{'}, u^{'})| \leq C_{m}\big(|z - z^{'}|+ |u - u^{'}|_{L^{2}(n)} \big). $$
Hence, for each $m$ and and $N$ being a fixed integer, we get existence of a solution  in $\mc{S}^{2} \times L^{2}(W) \times L^{2}(\tilde{N}_{p})$ of the BSDE given by ($f^{1,m}, \frac{B}{N}$): we denote it by ($Y^{1, m}, Z^{1, m}, U^{1, m}$).
% Referring to the a priori estimates in Proposition 2.2 \cite{Royer2}, then, for all $m$, $Y^{1, m}$ is in $S^{\infty}$ %and satisfies
%$$|Y_{s}^{1, m}|^{2}
% \leq \left|\frac{B}{N}\right|_{\infty}^{2} + \mb{E}\left(\int_{0}^{T} |f^{1,m}(s, 0, 0)|^{2}ds\right),
 %\quad \mb{P}\textrm{-a.s. and for all}\; s. $$
\item[2.] The sequence ($f^{1,m}$) is increasing and converges, $\mb{P}$-a.s and for all $s$, to $\tilde{f}$ in the following sense
$$ f^{1, m}(s, z, u) \nearrow  \tilde{f}(s, z, u), \quad \textrm{as} \; m \; \textrm{goes to} \; \infty. $$
%\item[3.] Each generator $f^{1,m}$ satisfies the same growth condition as $\tilde{f}$, i.e.
%\begin{equation}\label{eq: uniformbound}  - z\theta_{s} - \frac{|\theta_{s}|^{2}}{2\alpha} \le f^{1, m}(s, z, u) \le %\frac{\alpha}{2}|z|^{2} + |u|_{\alpha},
%\end{equation}
%and therefore, for each $m$, the triple ($Y^{1, m}, Z^{1, m}, U^{1, m}$) satisfies the a priori estimates stated in %lemma \ref{estim2}.
\end{enumerate}
Using both the Lipschitz property, the monotonicity of $(f^{1, m})$, the property ($H_{2}$) and the comparison result in Theorem 2.5 in \cite{Royer2}, ($ Y^{1, m}$) is increasing and hence, we can define $\tilde{Y}$ as follows
$$ \tilde{Y}_{s} := \disp{\lim \nearrow Y_{s}^{1, m} }, \quad \mb{P}\textrm{-a.s. and for all}\; s.$$
From the second assertion in Lemma \ref{estim2}, both the two sequences ($Z^{1, m}$) and ($U^{1, m}$) are bounded respectively in $L^{2}(W)$ and $L^{2}(\tilde{N}_{p})$: this entails the existence of weak limits denoted by $\tilde{Z}$ and $\tilde{U}$. \\

\indent To conclude this paragraph and for later use, we give a precise estimate of the norm of $Y^{1, m}$ in $\mc{S}^{\infty}$
\begin{equation}\label{eq: estimessent}
  |Y_{s}^{1, m}|_{\mc{S}^{\infty}} \leq \disp{ \frac{|B|_{\infty}}{N} },\quad \mb{P}\textrm{-a.s. and for all} \; s.
\end{equation}
(For sake of completeness, a detailed proof is provided in the first appendix A1.)
This estimate, which is independent of $m$, is essential in the proof of the monotone stability result given in the next paragraph: in particular, it allows to obtain the condition (\ref{eq: conditiontermF}) on $N$ under which the BSDE with parameters ($\tilde{f}, \; \frac{B}{N}$) admits a solution.
\\ 

\paragraph*{The stability result: convergence of the approximating sequence}
To justify that ($\tilde{Y}, \tilde{Z}, \tilde{U}$) solves the BSDE given by ($\tilde{f},\; \frac{B}{N}$), we prove the same kind of stability result as in \cite{mkobylanski} for the approximating sequence of BSDEs given by ($f^{1, m}, \frac{B}{N}$). To this end, we justify the three following convergence results  \\
(i) $\; \disp{Z^{1, m}  \to  \tilde{Z} \; (\textrm{in} \; L^{2}(W) ), \; \textrm{as} \; m \to \infty}$,\\
  (ii)$\; \disp{ U^{1, m}  \to  \tilde{U}, \;(\textrm{in} \; L^{2}(\tilde{N}_{p}(dx,ds))), \; \textrm{as} \; m \to \infty}$,\\
  (iii)  $\disp{\mb{E}\big(\int_{0}^{t}|f^{1, m}(s, Z_{s}^{1, m}, U_{s}^{1, m}) -\tilde{f}(s,\tilde{Z}_{s},\tilde{U}_{s}) | ds \big) \to 0, \; \textrm{as} \; m \to \infty}.$\\
Assertions (i) and (ii) correspond to the strong convergence of the sequences ($ Z^{1, m}$) and ($U^{1, m}$) to $\tilde{Z}$ and to $\tilde{U} $ in their respective Hilbert spaces. The proof being tedious and merely technical, it is relegated to the end in Appendix A2:
 we just give here the constraint condition on $N$: $M_{B}$ being an upper bound of $B$ in $L^{\infty}(\mc{F}_{T})$, $N$ should satisfy
 \begin{equation}\label{eq: conditiontermF} 
\frac{M_{B}}{N}  \le \disp{\inf \{ \frac{1}{32 \alpha}, \;\frac{1}{16C} \} },
\end{equation}
where $C$ is a constant depending only on $\alpha$ and $|B|_{\infty}$.\\
\indent To prove the convergence in $L^{1}(ds \otimes d\mb{P})$ stated in (iii), we apply the dominated convergence theorem by checking:\\
\begin{itemize}
\item $\;$ The convergence of $(f^{1, m}(s,Z_{s}^{1, m},U_{s}^{1, m}))$ to $\tilde{f}(s,\tilde{Z}_{s},\tilde{U}_{s})$, in $ds \otimes d\mb{P}$-measure, \\
\item $\;$ The existence of a uniformly integrable control of ($f^{1, m}(s,Z_{s}^{1, m},U_{s}^{1, m}))$  (independent of $m$). \\
\end{itemize}
The second assertion results easily from the inequality
\[  \begin{array}{l}
  |f^{m}(s,Z_{s}^{1, m}- \frac{\theta_{s}}{\alpha},U_{s}^{1, m})| \\
\\
\quad \quad \leq \disp{\max \left\{ \big(\frac{\alpha}{2}|Z_{s}^{1, m}- \frac{\theta_{s}}{\alpha}|^{2} + |U_{s}^{1, m}|_{\alpha}\big) ; \quad \big(-\theta_{s}(Z_{s}^{1, m}- \frac{\theta_{s}}{\alpha}) -\frac{|\theta_{s}|^{2}}{\alpha}\big)\right\}}.
\end{array} \]
To conclude for this second assertion, we rely on the uniform integrability of ($|Z^{1, m}- \frac{\theta}{\alpha} |^{2}$) and ($|U^{1, m}|_{\alpha} $), which results from their convergence in $L^{1}(ds \otimes d\mb{P})$ and on the boundedness assumption on $\theta$.
To prove the first point, we state an auxiliary result
\begin{proposition}\label{convuniform}
For all $s$ and for all converging sequences $(z^{m})_{m}$ and $(u^{m})_{m}$ respectively in $\mathbb{R}$ and $L^{2}( n(dx))$, such that the sequence ($u^{m}$) is uniformly bounded in $L^{\infty}(n)$ and satisfies:
$$\exists \; C > 0, \quad  \disp{\sup_{m}|u^{m}|_{L^{2}(n)}} \leq C,$$
 we have
$$ f^{1, m}(s, z^{m}, u^{m}) \to  \tilde{f}(s, z, u), \;  \mb{P}\textrm{-a.s. and for all}\; s, \; \textrm{as}\; \; m \to \infty . $$
\end{proposition}
The proof of this lemma results from the convergence of ($z^{m}$) and ($u^{m}$) (respectively to $z$ and $u$) and the simple convergence of $(f^{1, m})$ to $\tilde{f}$.\\
Without loss of generality and using the convergence results given in (i) and (ii), we can now assume\footnote{To ensure the convergence in $ds \otimes d\mb{P}$-measure, we ought to consider subsequences.} that both ($Z_{s}^{1, m}$) and ($U_{s}^{1, m} $) converge in $ds \otimes d\mb{P}$-measure to $\tilde{Z}_{s} $ and $\tilde{U}_{s} $ respectively in $\mb{R}$ and in $L^{2}(n)$: this entails the convergence in $L^{1}(ds \otimes d\mb{P})$ of ($f^{1, m}(s, Z_{s}^{1, m}, U_{s}^{1, m})$) to $\tilde{f}(s,\tilde{Z}_{s},\tilde{U}_{s}) $.\\
Passing to the limit in the equation satisfied by $Y^{1,m}$
\begin{equation}\label{eq: eqm}
 Y_{t}^{1, m} = \frac{B}{N} + \disp{\int_{t}^{T}f^{1, m}(s, Z_{s}^{1, m}, U_{s}^{1, m})ds - \int_{t}^{T}Z_{s}^{1, m}dW_{s} - \int_{t}^{T}\int_{\mathbb{R}\setminus\{0\}}U_{s}^{1, m}(x)\tilde{N}_{p}(ds, dx)}
\end{equation}
the increasing limit $\tilde{Y}$ satisfies
\begin{equation}\label{eq: eqlim}
 \tilde{Y}_{t} = \frac{B}{N} + \disp{\int_{t}^{T}\tilde{f}(s, \tilde{Z}_{s}, \tilde{U}_{s})ds - \int_{t}^{T}\tilde{Z}_{s}dW_{s} - \int_{t}^{T}\int_{\mathbb{R}\setminus\{0\}}\tilde{U}_{s}(x)\tilde{N}_{p}(ds, dx)}
\end{equation}
Substracting (\ref{eq: eqm}) and (\ref{eq: eqlim}) and taking then successively the supremum over $t$ and the expectation, we get
$$ \disp{\mb{E}\big(\sup_{t \in [0, T]} |Y_{t}^{1,m} - \tilde{Y}_{t}| \big) \to 0},$$
using, in particular, the Doob's inequalities for the square integrable martingales $ (Z^{1, m} - \tilde{Z}) \cdot W$ and $(\tilde{U} - U^{1, m}) \cdot \tilde{N}_{p}$ and the respective convergence of $ Z^{1, m} - \tilde{Z}$ in $ L^{2}(W) $ and $\tilde{U} - U^{1, m} $ in $ L^{2}(\tilde{N}_{p}(dx,ds))$.

\subsubsection{Step 2: the iterative procedure} 
In this step, we justify the existence result for the BSDE with parameters ($\tilde{f}, \; B$) ($B$ being an arbitrary bounded $\mc{F}_{T}$-measurable random variable).
%To this end, we explain the construction of an appropriate approximating sequence of BSDEs and we then prove %existence of solutions for the BSDEs introduced.\\ 
%(providing here a sufficient condition on the norm in $L^{\infty}$ of the terminal condition).\\

\paragraph*{Construction}
We provide here the explicit construction of a sequence of intermediate BSDEs with parameters ($f^{(i)} , \frac{B}{N} $) such as described at the beginning of section 3.3. For this, we define as follows the sequence $(f^{(i)})$: \\
\begin{enumerate}
 \item[1] We initialize by setting: $f^{(1)} := \tilde{f}$: the first step provides a solution for the BSDE with parameters ($f^{(1)}, \frac{B}{N}$) as soon as: $N \ge N^{1}$ with $N^{1}$ satisfying (\ref{eq: conditiontermF}). We denote this solution by ($\tilde{Y}^{1}, \tilde{Z}^{1}, \tilde{U}^{1}$). \\
 \item[2] Assuming that the sequence ($f^{(k)}$) is constructed up to step $k$, $k \ge 1$, and that each BSDE given by ($f^{(i)}, \frac{B}{N}$) (for an integer $N$ to give explicitely) admits a solution ($\tilde{Y}^{i}, \tilde{Z}^{i}, \tilde{U}^{i}$), we define the generator $f^{(k+ 1)}$ by setting
 $$f^{(k+ 1)}(s, z , u) = \tilde{f}(s,\;z + \bar{Z}_{s}^{k} - \frac{\theta_{s}}{\alpha}, \;u +\bar{U}_{s}^{k}) - \tilde{f}(s,\bar{Z}_{s}^{k}-\frac{\theta_{s}}{\alpha}, \bar{U}_{s}^{k}),$$ 
with: $\bar{Z}^{k} = \disp{\sum_{i \leq k}\tilde{Z}^{i}}$ and $\bar{U}^{k} = \disp{\sum_{i \leq k}\tilde{U}^{i}}$.\\
\end{enumerate}
 Provided there exists a solution ($\tilde{Y}^{i}, \tilde{Z}^{i}, \tilde{U}^{i} $) up to step $k$ and by definition of each $f^{(i)}$, we have: $ \disp{\sum_{i= 1}^{k}f^{(i)}(s, \tilde{Z}_{s}^{i}, \tilde{U}_{s}^{i}) = \tilde{f}(s, \bar{Z}^{k}, \bar{U}_{s}^{k})}$ and hence, the triple ($\bar{Y}^{k}, \bar{Z}^{k} ,\bar{U}^{k}$) with: $\bar{Y}^{k} = \disp{\sum_{i = 1}^{k} \tilde{Y}^{i}}$, solves the BSDE given by the generator $\tilde{f}$ and
the terminal condition equal to $\disp{\sum_{i= 1}^{k}\frac{B}{N}}$. After $N$ iterations of that procedure, it leads to a solution of the BSDE with parameters ($\tilde{f}, \; B$). \\
%It remains to prove the existence of a sequence of solutions ($\tilde{Y}^{i}, \tilde{Z}^{i}, \tilde{U}^{i} $)) (for %an appropriate choice of $N^{i}$) and a result of convergence analogous as in step 1.\\

\paragraph*{New stability result}
\paragraph*{Construction of the approximating sequence of BSDEs}
To justify the existence of a solution of the BSDE given by ($f^{(2)}, \;\frac{B}{N}$), we proceed analogously as in Section 3.3.1 by providing an explicit constraint on the integer $N$ (we deal with this technical issue in Appendix A3).
Keeping the same notation for $f^{m}$, we introduce \footnote{Assuming the procedure can be applied up to step $k$, then, for any $k$, $k \ge 2$, we define $f^{k+1, m}$ analogously
$$ f^{k+1, m}(s, z, u) := f^{m}(s, z + (\bar{Z}_{s}^{k,m} - \frac{\theta_{s}}{\alpha}), u +  \bar{U}_{s}^{k,m} ) - f^{m}(s,\bar{Z}_{s}^{k,m} - \frac{\theta_{s}}{\alpha}, \bar{U}_{s}^{k,m} ),$$
and since $(\bar{Z}^{k,m})$ (resp. $(\bar{U}^{k,m}) $) is uniformly bounded in $L^{2}(W) $ (resp. in $L^{2}(\tilde{N}_{p})$), the generator $f^{k+1, m}$ satisfies again the same growth condition and control of the increments as $f^{2,m}$. } the sequence $(f^{2, m})_{m}$ as follows
$$ f^{2, m}(s, z, u) := f^{m}(s, z + Z_{s}^{1,m} - \frac{\theta_{s}}{\alpha}, u +  U_{s}^{1,m}) - f^{m}(s,Z_{s}^{1,m} - \frac{\theta_{s}}{\alpha}, U_{s}^{1,m} ).$$
Using the same argumentation as in Step 1, we obtain a solution ($ Y^{2, m}, Z^{2, m}, U^{2, m}$)
of the BSDE given by ($f^{2,m}, \frac{B}{N}$). $f^{2, m}$ satisfying ($H_{1}$), both sequences ($ Z^{2, m}$) and ($ U^{2, m} $) are uniformly bounded respectively in $L^{2}(W) $ and in $L^{2}(\tilde{N}_{p})$ and we denote by $\tilde{Z}^{2}$ and $\tilde{U}^{2}$ their respective weak limits.\\
 By definition, the generator $f^{2, m}$ satisfies: $f^{2, m}(s, 0, 0) \equiv 0$ and hence, using the same procedure as described in Appendix A1, we get that any bounded solution $Y^{2, m}$ satisfies
 \begin{equation}\label{eq: newcontrol} |Y^{2, m}|_{\mc{S}^{\infty}} \le \big|\frac{B}{N}\big|_{\infty}.\end{equation}
%To prove this last estimate, we justify that $f^{2,m}$ satisfies an assumption of the same kind as ($H_{2}$) (since %it does not require any additional technicalities, we skip this and  
%We refer to the same procedure as in the first appendix to compute the increments of $f^{2, m}$ with respect to $z$ %and $u$.
 Now, to prove the existence of an almost sure limit for ($Y^{2,m}$), we cannot proceed as in step 1, since we do not have any monotonicity property for ($Y^{2, m}$): in fact, the sequence ($ f^{2, m}$) is neither increasing nor decreasing: however, if we consider $\bar{f}^{2,m}$ defined by: $\bar{f}^{2,m} = f^{2, m} + f^{1, m}$, then ($Y^{2, m} + Y^{1, m}$) is increasing and we can define $\bar{Y}^{2} $ as follows 
$$ \bar{Y}_{s}^{2}  = \disp{\lim_{m} \nearrow \big(Y_{s}^{2, m} + Y_{s}^{1, m}\big), \; \mb{P}\textrm{-a.s and for all} \; s. }$$
Since $(Y_{s}^{1, m})$ is increasing and converges to $\tilde{Y}_{s}$, $\mb{P}$-a.s. and for all $s$, $(Y_{s}^{2, m})$ converges to $\tilde{Y}_{s}^{2}$ defined by: $\tilde{Y}_{s}^{2}  = \bar{Y}_{s}^{2}- \tilde{Y}_{s}$.\\
The aim of the following paragraph is to prove a convergence result for the sequence ($ Y^{2, m}, Z^{2, m}, U^{2, m}$) and identify its limit ($\tilde{Y}^{2}, \tilde{Z}^{2}, \tilde{U}^{2}$) as a solution of the BSDE given by ($f^{(2)}, \frac{B}{N}$).\\
\paragraph*{Convergence of the approximating sequence}
As in Section 3.3.1, we have to prove the strong convergence of ($Z^{2, m} $) to $\tilde{Z}^{2} $ in $L^{2}(W)$ (respectively of ($U^{2, m}$) to $\tilde{U}^{2} $ in $L^{2}(\tilde{N}_{p})$)
and then justify a new stability result for the solutions of the BSDEs with parameters ($f^{2, m}, \; \frac{B}{N}$).\\
For sake of clarity, the proof of the strong convergence of ($Z^{2, m} $) and ($U^{2, m}$) is relegated to Appendix A3: using this last result and proceeding the same way as in the second paragraph in Section 3.3.1, we get 
$$\disp{\mb{E}\left(\sup_{t}|Y_{t}^{2, m} - \tilde{Y}_{t}^{2}| \right) +|Z^{2, m} - \tilde{Z}^{(2)}|_{L^{2}(W)} + |U^{2, m} - \tilde{U}^{(2)}|_{L^{2}(\tilde{N}_{p})}\to 0 } ,$$
and we identify the triplet ($\tilde{Y}^{(2)}, \tilde{Z}^{(2)}, \tilde{U}^{(2)}$) as a solution of the BSDE with parameters ($f^{(2)}, \frac{B}{N}$), $N$ satisfying (\ref{eq: nvcondittermF}) which is the new constraint\footnote{To obtain this constraint on the integer $N$, we rely on the fundamental estimate given by (\ref{eq: newcontrol}).} obtained in Appendix A3.

\paragraph*{End of the iteration procedure}
 In step 1, we have obtained a triple ($\tilde{Y}, \tilde{Z}, \tilde{U}$) solving the BSDE with parameters ($\tilde{f}, \frac{B}{N}$)
under the condition (\ref{eq: conditiontermF}) on $N$ and, in the previous paragraph, a solution ($ \tilde{Y}^{2}, \tilde{Z}^{2}, \tilde{U}^{2}$) of the BSDE with parameters ($f^{2}, \; \frac{B}{N}$) under the more restrictive condition (\ref{eq: nvcondittermF}). Defining $\bar{Y}^{2} $ by: $ \bar{Y}^{2} =\tilde{Y} + \tilde{Y}^{2}$ ($\bar{Z}^{2}$ and $\bar{U}^{2}$ being defined analogously), then ($\bar{Y}^{2}, \bar{Z}^{2}, \bar{U}^{2}  $) is solution of the BSDE given by ($\tilde{f}, \frac{2B}{N}$) (this holds if we choose for $N$ the minimal integer satisfying (\ref{eq: nvcondittermF})).\\
We distinguish two cases\\
\begin{itemize}
\item[1.] If we can choose $N = 2$, then the triple ($\bar{Y}^{2}, \bar{Z}^{2}, \bar{U}^{2}$) is the desired solution (of the BSDE with generator $\tilde{f}$ and terminal condition $B$).\\
\item[2.] In the second case, we proceed with at least one further iteration of the procedure described in step 2. For any $k$, $k \ge 2$, 
%and using the controls of the BMO norms of $\bar{Z}^{k}$ and $\bar{U}^{k}$ ,
we check that, for fixed $k$, each generator $f^{k, m}$, which is defined analogously as $f^{2, m}$ and whose expression is given at the bottom of page 14, satisfies an assumption similar to $(H_{2})$ and the property: $f^{k, m}(s, 0, 0) \equiv 0$. Under these two last assumptions, the following estimate holds for any $k$ and $m$
$$ |Y^{k, m}|_{\mc{S}^{\infty}} \le \frac{|B|_{\infty}}{N}.$$
Therefore, both the construction described in subsection 3.3.2 for the case $k = 2$ and the method to establish the stability result can be i-terated up to step $k$, $k \ge 2$ and
 in particular, at each step $i$, $i \ge 2 $, the condition (\ref{eq: nvcondittermF}) established in the second appendix remains unchanged. If we denote by $N^{1}$ the minimal integer satisfying (\ref{eq: nvcondittermF}) and if we then define ($Y, Z, U$) by: $(Y, Z, U) :=(\bar{Y}^{N^{1}}, \bar{Z}^{N^{1}}, \bar{U}^{N^{1}})$, with $\bar{Y}^{N^{1}} $ such that: $\bar{Y}^{N^{1}} = \disp{\sum_{i= 1}^{N^{1}}\tilde{Y}^{(i)}}$, this provides a solution of the BSDE with parameters ($\tilde{f}, B$).\\
%(hence, the  sequence  $(N^{i})_{i \ge 2}$ can be taken constant from step $2$). It entails the existence of $k$ such that $\disp{\sum_{i=1}^{k}\frac{1}{N^{i}} =1}$.
\end{itemize}

\subsubsection{Step 3: Conclusion}
In the previous steps, we have proved the existence of a solution of the BSDE (\ref{eq: eqavecsauts}) with parameters ($\tilde{f}$, $B$), where $B$ is an arbitrary bounded and
$\mc{F}_{T}$-measurable variable. Using this, we prove an existence result for
 the BSDE with parameters ($f$, $\bar{B}$), where the new terminal condition $\bar{B}$ can be expressed in terms of $B$.\\
Thanks to the two first steps, we can claim the existence of a triple ($Y, Z, U$) such that
\[ \begin{array}{ll}
   
Y_{t} = B + & \disp{\int_{t}^{T}[ f(s, Z_{s}-\frac{\theta_{s}}{\alpha},U_{s}) - f(s, -\frac{\theta_{s}}{\alpha}, 0)] ds  } \\
 & \; -\disp{ \;\int_{t}^{T}Z_{s}dW_{s} -\int_{t}^{T}\int_{\mathbb{R}\setminus\{0\}}U_{s}(x)\tilde{N}_{p}(ds, dx)},  \\
 \end{array} \]
which is well defined for any bounded random variable $B$.
If we define the processes $\bar{Y}$, $\bar{Z}$ and $\bar{U}$ as follows 
\begin{equation}\label{eq: expressionY}
\bar{Y}_{s} = \big(Y_{s} - \disp{\int_{0}^{s}f(u, -\frac{\theta_{u}}{\alpha}, 0)du -\int_{0}^{s}\frac{\theta_{u}}{\alpha}dW_{u}}\big),\; \; \bar{Z}_{s} = Z_{s} - \frac{\theta_{s}}{\alpha} \; \textrm{and} \; \bar{U}_{s} = U_{s},
\end{equation}
then, $\bar{Y}$ solves the following BSDE
$$\bar{Y}_{t} = \bar{B} + \disp{\int_{t}^{T}f(s, \bar{Z}_{s}, \bar{U}_{s})ds -  \int_{t}^{T}\bar{Z}_{s}dW_{s} -\int_{t}^{T}\int_{\mathbb{R}\setminus\{0\}}\bar{U}_{s}(x)\tilde{N}_{p}(ds, dx)},  $$
with generator equal to $f$ and terminal condition $\bar{B}$ equal to 
\begin{equation}\label{eq: expressionB} 
\bar{B} = B- \disp{\int_{0}^{T}f(s, -\frac{\theta_{s}}{\alpha}, 0)ds -\int_{0}^{T}\frac{\theta_{s}}{\alpha}dW_{s}}.
\end{equation}
 Due to (\ref{eq: expressionB}), the terminal condition $\bar{B}$ is no more in $L^{\infty}(\mc{F}_{T})$ and similarly, considering the first relation in (\ref{eq: expressionY}), $\bar{Y}$ is not in $\mc{S}^{\infty}$ but it only satisfies that $\exp(\bar{Y})$ is in $ \mc{S}^{p}$, for any $p$, $p > 0$.
To prove this, we use that
\begin{equation}\label{eq: relationcorrespondance} \exp(\alpha\bar{Y}_{t}) = \exp(\alpha Y_{t})\mc{E}\big(-\theta  \cdot W\big), 
\end{equation}
and we then rely on the boundedness of the process $\theta$ and on Novikov's criterion to obtain that $\mc{E}\big(- \theta \cdot W\big)$ admits moments of any order. 
 Since $Y$ is in $\mc{S}^{\infty}$, we obtain that $\bar{Y}$ admits exponential moments (the same holds for the terminal condition $\bar{B}$), which achieves the proof of (i) in Theorem 1.\\
Now, to obtain a solution for BSDE with parameters $f$ and $\bar{B}$, $\bar{B}$ being an arbitrary bounded random variable, we need to prove a more general existence result for BSDEs with generator $\tilde{f}$: this is the aim of the following section. 

\section{An existence result under more general condition}
In this section, we prove an existence result for solutions of BSDEs with generator $\tilde{f}$ and terminal condition $B$, under the restrictive condition that the terminal condition $B$ has exponential moments of any order: i.e.,
\begin{equation}\label{eq: existencemoments}
\forall \; \; k >0, \quad \quad \mb{E}\left(\exp(k |B|) \right) < \infty.
       \end{equation}
To prove a new existence result under this condition (\ref{eq: existencemoments}) on $B$, we adapt the procedure given in \cite{BriandetHu} in the context of a discontinuous setting and, for sake of clarity, we split the proof into three main steps.\\
Before proceeding with the proof, we give here the two properties ($H_{1}^{'}$) and ($ H_{2}^{'}$) satisfied by $\tilde{f}$. We first check that there exists a strictly positive constant $K$ and a non negative process $\bar{\alpha}$ satisfying: $\disp{\int_{0}^{T}\bar{\alpha}_{s}ds \le a} $, such that  
%$$ \exists\;  K >0, \; \bar{\alpha} \; \textrm{s.t.} \; \disp{\int_{0}^{T}\bar{\alpha}_{s}ds \le a}  $$
\begin{equation}\label{eq: growth}   (H_{1}^{'})\quad \quad   -\theta z \le \tilde{f}(s,\; z, \;u ) \le \bar{\alpha}_{s} + \frac{K}{2}|z|^{2} +|u|_{K}, 
 \nonumber
\end{equation}
which holds true when taking: $\bar{\alpha} = \frac{|\theta|^{2}}{\alpha}$ and $K = 2\alpha$. 
%(this also corresponds to the growth condition (which is analogous to ($H_{1}$))).\\ More precisely, for any $z^{1}, \;z^{2}$ in $\mb{R}$ and any $u^{1},\; u^{2}$ in $(L^{2} \cap L^{\infty})(n(dx))$, we have
Furthermore, the generator $\tilde{f}$ satisfies a new assumption denoted by ($H_{2}^{'}$) in the sequel and very similar to ($H_{2}$) stated in section 3.1 for the generator $f$. More precisely, for any $z^{1} $, $z^{2}$ in $\mb{R} $ and any ($u^{1} $, $u^{2}$) in $L^{2} \cap L^{\infty}(n)$, we have
\begin{itemize}
 \item[(1)] \[ \begin{array}{ll} \tilde{f}(s, z^{1}, u^{1}) - \tilde{f}(s, z^{2}, u^{1}) &= f(s, z^{1} -\frac{\theta_{s}}{\alpha}, u^{1}) - f(s, z^{2} -\frac{\theta_{s}}{\alpha}, u^{2})\\
\\
&  = \lambda^{'}(z^{1}, z^{2})(z^{1}-z^{2}),\\
              \end{array} \]
with $\lambda^{'}$ defined as follows 
\[\left\{  \begin{array}{ll}
     \lambda_{s}(z^{1}, z^{2}) \;& =  \frac{f(s, z^{1}, u) - f(s, z^{2}, u)}{z^{1} - z^{2}},  \; \textrm{if}\; z^{1} - z^{2} \neq 0,\\
     \\
    \lambda_{s}(z^{1}, z^{2}) &  = \; 0, \quad \textrm{otherwise}.\\
    \end{array} \right.
  \]  
and satisfying in particular that, as soon as $Z^{1}$ and $Z^{2}$ are in BMO($W$), the BMO property holds also for the process $ \lambda^{'}(Z^{1}, Z^{2})$.
%Such properties of $\lambda^{'}$ are justified in Appendix A1 (for the sequence of truncated generators %$(f^{1,m})$).
% and, using ($H_{2}$), it satisfies
%The process $ \lambda^{'}(z^{1}, z^{2})$ is equal to $\lambda(z^{1} -\frac{\theta}{\alpha}, %z^{2}-\frac{\theta}{\alpha})$: using assumption ($H_{2}$) satisfied by $f$ (we refer to section 3), we get 
%$$ \exists \; \kappa^{'} >0, \quad  |\lambda^{'}(z^{1}, z^{2})| \le C\big(\kappa^{'} + |z^{1}|+ |z^{2}|\big) \; \textrm{with}\; \kappa^{'} \le  \kappa + 2\frac{|\theta|}{\alpha}, $$
%$\kappa$ being introduced in ($H_{2}$).
\item[(2)] $$ \tilde{f}(s, z^{1}, u^{1}) - \tilde{f}(s, z^{1}, u^{2}) = \disp{\int_{\mb{R}^{*}}\gamma_{s}(u^{1}, \; u^{2})(u^{1} - u^{2})n(dx) },$$
where $\gamma$ has already been introduced in assumption ($H_{2}$) in Section 3.1.  
\end{itemize}

\paragraph*{Step 1: Comparison result and a priori estimates}
For later use, we provide here both a comparison theorem and a priori estimates. \\
\begin{proposition}\label{comparison}
Considering two bounded terminal conditions $\xi^{1}$ and $\xi^{2}$, if we denote by ($Y^{1}, Z^{1}, U^{1}$) (resp. ($ Y^{2}, Z^{2}, U^{2}$)) the solution in $\mc{S}^{\infty} \times L^{2}(W)\times L^{2}(\tilde{N}_{p})$ of the BSDE with parameters ($\tilde{f}, \;\xi^{1}$) (resp. ($\tilde{f}, \; \xi^{2}$)), then, as soon as: $ \xi^{1} \le \xi^{2}$, we have: $Y_{t}^{1} \le Y_{t}^{2}$, $\mb{P}$-a.s. and for all $t$.\\
\end{proposition}
Since the proof is based on the same ingredients as those given in Appendix A1, we skip the details and we just give the main steps:\\
$\bullet \;$ a standard linearization of the increments of the generator $\tilde{f}$  
$$\tilde{f}(s, Z_{s}^{1}, U_{s}^{1}) - \tilde{f}(s, Z_{s}^{2}, U_{s}^{2}), $$
obtained by relying on the assumption ($H_{2}^{'}$). \\
$\bullet \;$ an appropriate change of measure and a localization procedure to characterize $ Y^{1} -Y^{2}$ as a $\tilde{\mb{Q}}$-submartingale with terminal condition the non positive random variable $\xi^{1} - \xi^{2}$, for a suitable equivalent measure $ \tilde{\mb{Q}}$.\\ 
 
\begin{proposition}\label{aprioriestimates}
If we consider a BSDE with generator satisfying ($H_{1}^{'}$) and bounded terminal condition $B$, then, for any solution in $\mc{S}^{\infty} \times L^{2}(W)\times L^{2}(\tilde{N}_{p}) $, we have
\begin{equation}\label{eq: aprioriestimates}
\exists \;a,\; K > 0, \;C \; \textrm{s.t.}\; \; -C \mb{E}\big(|B|^{2}|\mc{F}_{t} \big)^{\frac{1}{2}} \le \bar{Y}_{t} \le \frac{1}{K}\ln\mb{E} \left(\exp(K (B + a))|\mc{F}_{t}\right),\end{equation} 
where the constant $K$ can be taken equal to 2$\alpha$, the constant $C$ can be taken equal to the norm in $\mc{S}^{2}$ of the stochastic exponential $\mc{E}(- \theta \cdot W)$ \footnote{To justify that the stochastic exponential $\mc{E}(- \theta \cdot W)$ is in $\mc{S}^{2}$, we use Novikov's criterion.} and the constant $a$ already introduced in ($H_{1}^{'}$) corresponds to an upper bound of $\disp{\int_{0}^{T}\bar{\alpha}_{s}ds} $.\\
\end{proposition}
Since it is very similar as in \cite{morlais2}, we only give the main ingredients: for the upper bound, it relies both on the application of It\^o's formula to $\exp(K Y)$ and on standard computations.
For the estimate in the left-hand side, we use that the lower bound of $\tilde{f}$ has linear growth with respect to its variable $z$ and that $\tilde{f}$ is such that: $\tilde{f}(s, 0, 0) \equiv 0$. 
Hence, $Y$ is greater than the solution of the linear BSDE with generator $-\theta z$ and terminal condition $B$, which is equal to $ \mb{E}^{\mb{P}^{\theta}}(B| \mc{F}_{t})$, with $\frac{d\mb{P}^{\theta}}{d\mb{P}} = \mc{E}(-\theta \cdot W)$: the terminal condition being bounded (and hence square integrable), a lower estimate is given by the expression in the left-hand side in (\ref{eq: aprioriestimates}).\\
\begin{flushright}
 $\square$
\end{flushright}

\paragraph*{Step 2: the stability result}
In this paragraph, we explain the construction of a sequence of BSDEs and for this sequence, we prove an extended stability result. For this, we make use of a localization procedure which is analogous as in \cite{BriandetHu}. \\
Our first aim is to obtain uniform a priori estimates, for any sequence of solutions ($\bar{Y}^{n}, \bar{Z}^{n}, \bar{U}^{n}$) of BSDEs with parameters ($\tilde{f}, B^{n}$), when the sequence $(B^{n})$ of terminal conditions is uniformly bounded in $\mc{S}^{\infty}$. 
\\
Assuming that $B$ is non negative\footnote{For the general case, we refer to \cite{BriandetHu}: setting first: $B^{n,\;p} =B \wedge n - (-B \wedge p)$, we construct a sequence $(\bar{Y}^{n,\; p})$ of solutions of the BSDEs given by ($\tilde{f}, B^{n, \;p}$) such that it is decreasing w.r.t $p$. The next step consists in establishing a stability result for this decreasing sequence, which is skipped here since it is analogous to the proof of Lemma \ref{newstability} and relies on the same kind of localization procedure and on the lower estimate obtained in (\ref{eq: aprioriestimates}).} and satisfies (\ref{eq: existencemoments}), we first define ($B^{n}$) as follows: $B^{n} = B \wedge n$. Using the results of Section 3, the BSDE with parameters $\tilde{f}$ and $B^{n}$ has a solution ($\bar{Y}^{n}, \bar{Z}^{n}, \bar{U}^{n}$) such that $\bar{Y}^{n}$ is in $\mc{S}^{\infty}$.
Thanks to the priori estimates given by (\ref{eq: aprioriestimates}) in Lemma \ref{aprioriestimates} and using that $B^{n}$ satisfies: $0 \le B^{n} \le B $, we obtain
$$ 0  \le \bar{Y}_{t}^{n}  \le  \frac{1}{K}\ln\mb{E} \left(\exp\big(K (B + a)\big)|\mc{F}_{t}\right),$$
where the expression of $K$ is explicited in the first step. Due to assumption (\ref{eq: existencemoments}), the random variable in the right-hand side is almost surely finite.
% and this property holds for the random variable in the left-hand side (for this, we argue that any random variable %having a finite positive exponential moment is square integrable).\\
\indent The first step of the localization procedure consists in introducing a sequence ($\tau_{k}$)$_{k}$ of stopping times as follows 
$$\tau_{k} = \disp{\inf\{t, \;\frac{1}{K}\ln\mb{E} \left(\exp(K (B+a))|\mc{F}_{t}\right) \ge k  \} \wedge T.} $$
If we then fix $k$ and if we denote by $\bar{Y}^{k, n}$ the process $\bar{Y}^{n}$ stopped at time $\tau_{k}$, this process solves a BSDE with generator $\tilde{f}^{k} = \tilde{f}\mathbf{1}_{\tau^{k} \le T}$ and terminal condition $\xi^{n, k}$ defined by
\[ \xi^{n, k}  = \left\{ \begin{array}{l}
      B^{n}, \; \textrm{if} \; \tau_{k} = T,\\
\\
     \bar{Y}_{\tau_{k}}^{n},\; \textrm{if} \; \tau_{k} < T. \\   
   \end{array} \right. \]

We now state a new stability result\footnote{For a very similar result in the brownian setting, we also refer to Lemma 3 in \cite{BriandetHu}.} for the sequence ($\bar{Y}^{k, n},\;\bar{Z}^{k, n} ,\;\bar{U}^{k, n}$) of solutions of the BSDEs with parameters ($\tilde{f}^{k}, \xi^{n, k} $), $k$ being fixed.
%We first state the stability result 
 
\begin{proposition}\label{newstability}
 Under the two following assumptions on the sequence of BSDEs with parameters ($f^{n} ,\; \xi^{n, k}$)$_{n}$\\
$ \bullet \;$ for all $n$, $f^{n} = \tilde{f}^{k}$, with $\tilde{f}^{k}$ satisfying assumption ($H_{1}^{'}$),\\
$ \bullet \;$ $(\xi^{n, k})$ is increasing and uniformly bounded in $\mc{S}^{\infty}$,\\
and if, in addition, there exists a sequence ($\bar{Y}^{k, n},\;\bar{Z}^{k, n} ,\;\bar{U}^{k, n}$) of solutions for the BSDEs with parameters ($\tilde{f}^{k},  \; \xi^{n,k}$) such that $(\bar{Y}^{k, n})$ is increasing then, there exists a triple ($\bar{Y}^{k}, \bar{Z}^{k}, \bar{U}^{k} $) such that
\begin{equation}\label{eq: convergenceprocessus}
 \mb{E}\left(\disp{\sup_{[0, T]} |\bar{Y}_{t}^{k, n} - \bar{Y}_{t}^{k}|}  \right) + \;|\bar{Z}^{n,k} -\bar{Z}^{k}|_{L^{2}(W)}  + \;|\bar{Z}^{n,k} -\bar{U}^{k}|_{L^{2}(\tilde{N}_{p})} \to 0,\end{equation} 
and this triple
solves the BSDE given by ($\tilde{f}^{k}, \xi^{k}$) (with $\xi^{k} $ defined by: $ \xi^{k} = \disp{\sup_{n}\xi^{n, k}} $).\\
\end{proposition}

 To justify the stability result stated in lemma \ref{newstability}, we apply the same procedure as in Appendix A2.
We first check all the required assumptions: by definition, $(\xi^{n, k})$ is an increasing sequence of bounded terminal conditions such that: $\disp{\sup_{n}|\xi^{n, k}|} \le k$ and, for all $n$, the generator $f^{n}$ equal to $ \tilde{f}^{k}$ satisfies the same assumptions than $\tilde{f}$: hence, we deduce \\
$\bullet \;$ the sequence
  $(\bar{Y}^{k, n})$ is increasing (this results from the comparison result given in lemma \ref{comparison}),\\
$\bullet \;$ $(\bar{Y}^{k, n})$ is uniformly bounded in $\mc{S}^{\infty} $ (i.e., the bounds are independent of $n$) with
$$\disp{ 0 \le  \sup_{n}  \bar{Y}^{k, n}\le k }. $$
Hence, we can define the process $\bar{Y}^{k}$ as follows\\
$$ \bar{Y}^{k} = \disp{\lim \nearrow_{k} \bar{Y}^{k, n}}.$$
Using standard computations (which are similar as in the proof of Lemma 3), we obtain that the two sequences $(\bar{Z}^{k, n})$ and $(\bar{U}^{k, n}) $ are bounded respectively in $L^{2}(W)$ and in $L^{2}(\tilde{N}_{p})$ and we denote by $\bar{Z}^{k} $ and by $\bar{U}^{k}$ their respective weak limits.\\
%Since ($\bar{Y}^{k, n}$) is uniformly bounded and using both the It\^o formula and standard computations (similar as %in the proof of (ii) in Lemma \ref{estim2}), we get that the sequences $(\bar{Z}^{k, n} )$ and $(\bar{U}^{k, n}) $ %are bounded respectively in $L^{2}(W)$ and in $L^{2}(\tilde{N}_{p})$: we denote by $\bar{Z}^{k} $ and by %$\bar{U}^{k}$ their respective weak limits.\\
To prove the strong convergence of both ($\bar{Z}^{k, n} $) and ($\bar{U}^{k, n} $), we follow the same procedure as in Appendix A2. For this, we apply the It\^o's formula to $ |\bar{Y}_{\cdot}^{k, n} - \bar{Y}_{\cdot}^{k, m}|^{2}$ and we rely on the following estimate
$$ |\bar{Y}_{\cdot}^{k, n} - \bar{Y}_{\cdot}^{k, m}|_{\mc{S}^{\infty}} \le |\xi^{n} - \xi^{m}|_{\infty}.$$
To justify this claim, we proceed as in Appendix A1: for this, we use that, for any $k$, the generator $\tilde{f}^{k}$ satisfies the same kind of assumption as $\tilde{f}$, that is ($H_{2}^{'}$) and we follow the same method as described in Appendix A1 to prove that $ \bar{Y}_{\cdot}^{k, n} - \bar{Y}_{\cdot}^{k, m} $ is a bounded $\mb{Q}$-submartingale with terminal condition equal to $\xi^{n} - \xi^{m}$ (for a well chosen equivalent measure $\mb{Q}$).\\ 
As a consequence, to rewrite the proof given in Appendix A2, we only need to check the sufficient condition  
\begin{equation}\label{eq: contrainte}
\exists \; M, \quad \disp{\sup_{n, m \ge M} |\xi^{n} - \xi^{m}|_{\mc{S}^{\infty}}} \le \disp{\inf \{\frac{1}{32 \alpha}, \;\frac{1}{16C} \}}.
\end{equation}
(This condition is obtained for a constant $C$ depending only on the parameters of the BSDE).
Since ($\xi^{n}$) converges in $L^{\infty}(\mc{F}_{T})$, it is a Cauchy sequence and, provided we take $M$ large enough, condition (\ref{eq: contrainte}) is ensured.
Hence, there exists a triple ($\bar{Y}^{k}, \bar{Z}^{k}, \bar{U}^{k}$) such that (\ref{eq: convergenceprocessus}) holds and solving
the BSDE with parameters $\tilde{f}^{k}$ and terminal condition $\xi^{k}$ such that:
$ \xi^{k} = \disp{\sup_{n} \xi^{n, k} }. $

\paragraph*{Step 3: conclusion}
We first define $Y$, $Z$ and $U$ as follows
$$Y_{t} = \bar{Y}_{t}^{k}\mathbf{1}_{t \le \tau^{k}}, \;Z_{t} = \bar{Z}_{t}^{k}\mathbf{1}_{t \le \tau^{k}} \; \textrm{and}\; U_{t} = \bar{U}_{t}^{k}\mathbf{1}_{t \le \tau^{k}}. $$
and to ensure the consistency of this definition, we need to check 
\begin{equation}\label{eq: consistency}
 \bar{Y}^{ k} \equiv \bar{Y}^{ k+1}  \quad \quad \textrm{on} \; [0, \tau^{k}]. 
\end{equation}
For this, we claim that, for each $n$ and each $k$, the solution ($\bar{Y}^{n, k}, \bar{Z}^{n,k}, \bar{U}^{n,k}$) of the BSDE with parameters ($\tilde{f}^{k}, \;B^{n}$) is unique\footnote{This uniqueness result follows from the comparison result stated in Lemma \ref{comparison}.}. Using then that the generators $\tilde{f}^{k} $ and $ \tilde{f}^{k + 1}$ coincide on $[0, \tau^{k}]$, we necessarily have: $\bar{Y}^{ k, n} \equiv \bar{Y}^{ k+1, n} $ on $ [0, \tau^{k}]$ and (\ref{eq: consistency}) results from the fact that $\bar{Y}^{ k} $ and $\bar{Y}^{ k+1} $ are the increasing and almost sure limits of ($\bar{Y}^{n, k} $) and ($ \bar{Y}^{n, k+1}$).\\
 Furthermore, since $B$ satisfies the property given by (\ref{eq: existencemoments}), the sequence ($\tau^{k}$) is stationnary (almost surely): this means that, for almost $\omega$, there exists $k(\omega) $ such that $\tau^{k}(\omega) = T$ and hence: $\xi^{k(\omega)} =B$.
As a consequence, the triple ($Y, Z, U$) solves the BSDE with parameters ($\tilde{f}, B$).\\

\indent To conclude, we use the result of the previous section: i.e, the existence of solutions of the BSDE with parameters $f$ and $\bar{B}$, for any random variable $\bar{B}$ defined in terms of $B$ as follows
\begin{equation}\label{eq: correspondance}
 \bar{B} = B- \disp{\int_{0}^{T}f(s, -\frac{\theta_{s}}{\alpha}, 0)ds -\int_{0}^{T}\frac{\theta_{s}}{\alpha}dW_{s}}.\end{equation}
(this expression is given in the last step in section 3). But in general, when $B$ is bounded, the random variable $\bar{B}$ is no more bounded (it only satisfies (\ref{eq: existencemoments})). 
To obtain the desired existence result, we consider an arbitrary bounded random variable $\bar{B}$ and we define $B$ in terms of $\bar{B}$ using (\ref{eq: correspondance}). Such a random variable $B$ satisfies the property (\ref{eq: existencemoments}) and hence, using the new existence result proved in this section, we obtain a solution ($Y, Z, U$) of the BSDE with parameters ($\tilde{f}, B$).
%To obtain a solution for the BSDE given by ($f, \bar{B}$) under the assumption of boundedness on $\bar{B}$, 
%(\ref{eq: correspondance}) shows that we need to consider an auxiliary BSDE with parameters ($\tilde{f}, B$) such %that $B$ has only exponential moments of any order. Using the new existence result of this section, we get a solution ($Y, Z, U$) of the BSDE with parameters ($\tilde{f}, B$). 
Defining then $(\bar{Y}, \bar{Z}, \bar{U})$ as follows
$$ \bar{Y}_{s} = \big(Y_{s} - \disp{\int_{0}^{s}f(u, -\frac{\theta_{u}}{\alpha}, 0)du -\int_{0}^{s}\frac{\theta_{u}}{\alpha}dW_{u}}\big),\; \; \bar{Z}_{s} = Z_{s} - \frac{\theta_{s}}{\alpha} \; \textrm{and} \; \bar{U}_{s} = U_{s},
$$ 
 this triplet solves the BSDE with parameters ($f$, $\bar{B}$). Since $\bar{B}$ is a bounded random variable and since $f$ satisfies $(H_{1})$, the application of Lemma \ref{estim2} entails that $\bar{Y}$ is in $\mc{S}^{\infty}$, which achieves the proof of (ii) in Theorem \ref{existence}.
%Conclude for existence of a solution for BSDE given by ($f, \; B$) for any bounded terminal condition.
%TO WRITE
\section{Application to the utility maximization pro-blem}
In this section, we make use of the notations introduced in Section 2 and using the results of the two previous sections, we provide a characterization of the value process at time 0
\begin{equation}\label{eq: pboptim}
\nonumber  V(x) = \displaystyle{\sup_{\pi \in \mathcal{A}}{\mathbb{E}(U_{\alpha}(X_{T}^{\pi} - \bar{B}))}},
\end{equation}
which is associated to the classical utility maximization problem with bounded liability $\bar{B}$.
% this random variable is given in terms of both the bounded random variable $B$ and the parameters $ \theta$ and $ %\alpha$ by (\ref{eq: expressionB}). In fact, using the study of the previous section and to solve the problem by %adapting arguments in \cite{morlais2}, we ought to consider this specific form of liability. 
We now
state the main result of this section.
\begin{theorem}\label{expressionprocessV}
The expression of the value process at time 0 is given by 
\begin{equation}\label{eq: fonctionV} \disp{ V(x) = -\exp(- \alpha(x - \bar{Y}_{0})),} \end{equation}
where $\bar{Y}_{0}$ represents the initial data of the solution ($\bar{Y}, \bar{Z}, \bar{U}$) to the BSDE (\ref{eq: eqavecsauts}) given by the parameters ($f$, $\bar{B}$) with the generator $f$ defined as follows
\begin{eqnarray*}
 \nonumber \disp{ f(s,z,u) = \disp{\inf_{\pi \in
 \C}\left(\frac{\alpha}{2}|\pi\sigma_{s} - (z+ \frac{\theta}{\alpha})|^{2} + |u - \pi \beta_{s}|_{\alpha} \right) -
  \theta z - \frac{|\theta|^{2}}{2\alpha}}  }. \\
    \nonumber
\end{eqnarray*}
Moreover, there exists an optimal and admissible strategy $\pi^{*}$, such that: $\pi^{*} \in \mc{A}$. Such a strategy satisfies $\mb{E}(U_{\alpha}(X_{T}^{\pi^{*}} - \bar{B})) = V(x),$ and it is characterized by
\begin{eqnarray}\label{eq: optimstrat}
     \disp{\pi^{*}_{s}(\omega) \in  \textrm{arg}\displaystyle{\min_{\pi \in \C}\big(\frac{\alpha}{2}|\pi\sigma_{s} - (Z_{s}+ \frac{\theta_{s}}{\alpha})|^{2} + |U_{s} - \pi \beta_{s}|_{\alpha}\big) }, \; \; \mathbb{P}\textrm{-a.s.}  \;\textrm{and} \; \textrm{for all} \; s }\\
  \nonumber   \end{eqnarray}
\end{theorem} 
%In that theorem and using that $\bar{Y}$ is in $\mc{S}^{p}$, for any positive $p$, this implies that $\exp(\alpha %\bar{Y})$ is at least in $L^{1}(ds \otimes d\mb{P})$ and hence, the expression (\ref{eq: fonctionV}) for the value %function is almost surely finite.
 Since it relies on the same procedure as in \cite{morlais2}, we give here a brief proof with the main arguments.\\

\subsubsection*{Proof of theorem \ref{expressionprocessV}} 
\hspace{0.5cm} We first denote by ($\bar{Y}, \bar{Z}, \bar{U}$) the solution in $\mc{S}^{\infty} \times L^{2}(W) \times L^{2}(\tilde{N}_{p}) $ of the BSDE given by ($f$, $\bar{B}$) whose existence has been obtained in the previous sections and, for any admissible $\pi$, we define $R^{\pi}$ as follows  
\begin{equation}\label{eq: expressionR}
  \forall \; t, \; \;  R_{t}^{\pi} = - e^{-\alpha X_{t}^{\pi}}e^{\alpha \bar{Y}_{t}}.
\end{equation}
 In a first step and to obtain the expression (\ref{eq: fonctionV}), 
 we prove the supermartingale property of $R^{\pi}$, which holds for any admissible strategy $\pi$ ($\pi \in \mathcal{A}$). 
Using standard computations derived from the It\^o's formula, $R^{\pi}$ has the following product form
$$ R_{t}^{\pi} =  R_{0}^{\pi}\tilde{M}_{t}^{\pi}e^{A_{t}^{\pi}}, $$
with the process $\tilde{M}$ such that
$$ \tilde{M}_{t} = \mathcal{E}_{t}(M) = \disp{\mathcal{E}_{t} \left((-\alpha (\pi \sigma - Z) \cdot W) + (e^{(- \alpha(\pi \beta -
  U))} - 1) \cdot \tilde{N}_{p})\right)},$$ 
and with $M^{\pi}$ and $A^{\pi}$ defined by: $ M^{\pi} = (-\alpha (\pi \sigma - Z) \cdot W) + (e^{(- \alpha(\pi \beta -
  U))} - 1) \cdot \tilde{N}_{p}, $
and by
$$A_{t}^{\pi} = \int_{0}^{t} \alpha \left (-\pi_{s} b_{s} - f(s,Z_{s},U_{s}) + \frac{\alpha}{2}|\pi_{s} \sigma_{s} - Z_{s}|^{2} + |U_{s} -\pi_{s} \beta_{s}|_{\alpha} \right)ds. $$
\indent Since $\tilde{M}^{\pi}$ is a non negative stochastic exponential, it is a local martingale for any $\pi$, and
 consequently, there exists a sequence of stopping times $(\tau^{n})$ converging to $T$ such that $\tilde{M}_{. \wedge\tau^{n}}^{\pi}$ is a martingale. By definition of the generator $f$, $\exp(A^{\pi})$ is non decreasing and since $R_{0}$ is non positive, $R_{\cdot \wedge\tau^{n}}^{\pi }$ satisfies
 \begin{equation}\label{eq: supermproperty} 
\forall A \in \mc{F}_{s}, \quad \mathbb{E}(R_{t \wedge\tau^{n}}^{\pi}\mathbf{1}_{A}) \leq \mathbb{E}(R_{s \wedge\tau^{n}}^{\pi}\mathbf{1}_{A}),
 \end{equation}
Using the definition (\ref{eq: expressionR}) of $R^{\pi}$, the uniform integrability of ($R_{. \wedge\tau^{n}}^{\pi}$)$_{n}$ results both from the uniform integrability of $e^{- \alpha X^{\pi}}$ (proved in Lemma \ref{classequality}) and the boundedness of $\bar{Y}$. Hence, passing to the limit as $n$ goes to $\infty$ in (\ref{eq: supermproperty}), it implies that, for all $A \in \mathcal{F}_{s}$, $\mathbb{E}(R_{t}^{\pi}\mathbf{1}_{A}) \leq \mathbb{E}(R_{s}^{\pi}\mathbf{1}_{A})$, which yields the supermartingale property of $R^{\pi}$.\\

%To prove the desired supermartingale property, we show the integrability in $L^{1}(ds \otimes d\mb{P})$ of the %process ($R^{\pi}$) (for any admissible strategy $\pi$): using that: $R^{\pi} = \exp(-\alpha %X^{\pi})\exp(\alpha\bar{Y})$, we get
%$$ \mb{E}\left( \big| \exp(-\alpha X_{t}^{\pi})\exp(\alpha\bar{Y}_{t})\big| \right) \le \mb{E}\left( \exp(|p\alpha %X_{t}^{\pi}|) \right)^{\frac{1}{p}} \mb{E}\left(\exp(q\alpha\bar{Y}_{t}) \right)^{\frac{1}{q}}< \infty,$$
%where, in that previous expression, the constant $p$ ($p > 1$) is obtained using Lemma \ref{classequality} and $q$ is %uniquely defined such that: $\frac{1}{p} + \frac{1}{q} = 1$. Since $\bar{Y}$ is in $\mc{S}^{q}$ (for any $q$), the %conclusion follows. We then deduce the uniform integrability of ($R_{\cdot \wedge\tau^{n}}^{\pi } $) and taking the %limit as $n$ goes to $\infty$ in (\ref{eq: supermproperty}), we obtain the supermartingale property of ($R^{\pi}$).\\

\indent To complete the proof of this theorem and justify the expression (\ref{eq: fonctionV}) for $V$, we first prove the optimality of any strategy $\pi^{*}$ satisfying (\ref{eq: optimstrat}).
From this last characterization of $\pi^{*}$, we obtain: $A^{\pi^{*}} \equiv 0$ and this entails that $R^{\pi^{*}} $ such that: $R^{\pi^{*}} = R_{0}^{\pi^{*}}\tilde{M}^{\pi^{*}}$, is a local martingale. By its definition, $\pi^{*}$ takes its value in $\mc{C}$ and hence, thanks to Lemma 1, $\pi^{*}$ is in $\mc{A}$, which entails that $R^{\pi^{*}}$ is a true martingale. From this last martingale property, we get
$$ \disp{\sup_{\pi \in \mc{C}} \mb{E}(R_{T}^{\pi})} = \mb{E}(R_{T}^{\pi^{*}}) = R_{0} = -\exp\big(-\alpha(x - \bar{Y}_{0})\big), $$
which gives the expression (\ref{eq: fonctionV}) for $V$.\\ 
\begin{flushright}
$ \square$
\end{flushright}

\section{Conclusion}
In this paper, we consider the utility maximization problem with an additional liability and under portfolio constraints in the context of a discontinuous filtration: we then solve this problem by using the same methodology than in \cite{ImkelleretHu}: this consists in relying both on the dynamic programming principle and on BSDEs techniques to obtain the expression of the value process in terms of the solution of a quadratic BSDE with jumps. Furthermore, since we relax the finiteness assumption of the Levy measure, this study extends of the results already obtained in \cite{morlais2}: under this additional restriction, we establish a new existence result, which is the main achievement of this paper. Then and as in \cite{morlais2}, this theoretical study allows to characterize explicitely and dynamically the value process associated to the utility maximization problem and also to 
prove existence of optimal strategies.
%in particular, we show that existence for the BSDE already introduced in \cite{morlais2} can be obtained provided we %restrict the class of terminal condition: as a consequence, an expression of the value process of the utility %maximization problem is given for specific form of liabilities.

\newpage
\section{Appendix}
\subsection{A1: Proof of the estimates (\ref{eq: estimessent}) and (\ref{eq: newcontrol})} 
Our aim here is to justify that, for any solution of the BSDE with parameters ($f^{k, m}, \; \frac{B}{N}$), we have
$$ |Y_{s}^{k, m}|_{\mc{S}^{\infty}} \leq \disp{ \frac{|B|_{\infty}}{N} },\quad \mb{P}\textrm{-a.s. and for all} \; s.$$
The cases when $k=1$ and $k=2$ corresponds to the inequalities (\ref{eq: estimessent}) and (\ref{eq: newcontrol}) already stated in section 3.3.1 and 3.3.2 and of great use in the proof of the two stability results in Appendix A2 and A3.\\
% these stability results are associated to the sequences ($Y^{k, m}, Z^{k, m}, U^{k, m}$).\\ 
\indent In this paragraph, we only consider the case when $k= 1$ (in fact, the general case is based on the same procedure, provided we check that for all $k$ and $m$, the increments of the generator $f^{k,m}$ satisfy analogous controls as those which are stated in ($H_{2}$) in Section 3.1 for $f$ or in Section 4 for $\tilde{f}$.)\\
\indent Now and in a first step, we proceed with the proof of the upper bound for $Y^{1,m}$ and, for this, we make use of a standard linearization procedure which we are going to describe.
%this consists in checking that, for each $m$, $f^{1,m}$ satisfies an assumption which is very similar to ($H_{2}$) %and then in using these controls of the increments to linearize the expression appearing in the It\^o's formula.\\
Firstly, for any $z, z^{'}$ in $\mb{R}$, $u, \; u^{'}$ in $\big(L^{2} \cap L^{\infty}\big)\big(n\big)$, we check
$$ f^{1, m}(s, z, u) - f^{1, m}(s, z^{'}, u^{'}) =f^{m}(s, z- \frac{\theta}{\alpha}, u) - f^{m}(s,z^{'}- \frac{\theta}{\alpha}, u^{'}) , $$
and therefore, we only need to consider the increments of the function $\tilde{f}^{m}$ defined by: $\tilde{f}^{m}:(s, z, u) \to f^{m}(s, z- \frac{\theta}{\alpha}, u)$.
Concerning the increments w.r.t. $u$, the upper bound given in ($H_{2}$) in Section 3.1 holds again (with the same process $\gamma$). For the increments w.r.t. $z$, we rewrite $f^{m}$ as follows 
$$ f^{m}(s, z, u) = \disp{\inf_{\pi \in \mc{C}} \left(\Phi(z, \pi)\rho_{m}(z) + \disp{\int_{\mb{R}^{*}}g_{\alpha}^{m}(u - \pi \beta) n(dx)}\right)}, $$
 with the function $\Phi$ which is defined by: $\Phi(z) = \Phi(z, \pi) =\frac{\alpha}{2}|\pi \sigma - (z+\frac{\theta}{\alpha}) |^{2}$ and is a continuously differentiable function whose differential has linear growth w.r.t. $z$. We also rely on
$$\disp{\inf_{\pi \in \mc{C}}F(\pi, z, u)} -\disp{\inf_{\pi \in \mc{C}}F(\pi, z^{'}, u)} \le \disp{\sup_{\pi \in \mc{C} }|F(\pi, z, u) - F(\pi, z^{'}, u) |}, $$
 and we then use an explicit upper bound for the increments of: $z \to  \Phi(z)\rho_{m}(z)$ to obtain 
\begin{equation}\label{eq: increments} |f^{m}(s, z- \frac{\theta_{s}}{\alpha}, u) -f^{m}(s, z^{'}- \frac{\theta}{\alpha}, u)| \le \end{equation}
$$\quad \quad \quad  \disp{|\sup_{\pi \in \mc{C}}{\left(\disp{\sup_{\lambda \in [0, 1]}} \Phi^{'}(z^{\lambda})\rho^{m}(z^{\lambda}) + \Phi(z^{\lambda})\big(\rho^{m}\big)^{'}(z^{\lambda})\right)}| z- z^{'}|  }, $$
$\left(\; \textrm{with:}\; z^{\lambda} = \lambda(z- \frac{\theta_{s}}{\alpha}) + (1 -\lambda)\big(z^{'}- \frac{\theta}{\alpha}\big)\right)$.
Using then that $(\rho^{m})^{'}$ is equal to zero except on $[m, m+1]$ (where it is bounded since continuous) and the increasing property of $\Phi$ (on $[m, m+1]$), we get $$\exists \; C> 0, \; \forall \; \lambda \in [0, 1], \quad  |\Phi(z^{\lambda})\big(\rho^{m}\big)^{'}(z^{\lambda})| \le C\Phi(m+1).$$
 Due to the assumptions on the parameters and the compactness of $\mc{C}$, the term in the right-hand side is a bounded process (we denote by $C_{m}$ an upper bound).
 Relying now on the linear growth of $\Phi^{'}$, straightforward computations leads to
\begin{tabbing}
$ $\= $ |f^{m}(s, z- \frac{\theta_{s}}{\alpha}, u) -f^{m}(s, z^{'}- \frac{\theta}{\alpha}, u)|$\\
$ $\>  $\quad \quad \quad \le \left(C_{m} + \; \disp{\sup_{\lambda \in [0, 1]}\Phi^{'}(z^{\lambda})\rho^{m}(z^{\lambda})} \right)|z- z^{'}|,$\\
$ $\>  $\quad \quad \quad \le C\left( \kappa^{m} + |z| + |z^{'}|\right)|z- z^{'}|,$\\ \end{tabbing}
with $\kappa^{m}$ in BMO($W$) and depending only on the parameters $\alpha$, $\theta$ and on $m$.\\
Defining $\lambda^{m}$ the same way as in Section 3.1 as follows
\[\left\{  \begin{array}{ll}
     \lambda^{m}_{s}(z, z^{'}) \;& :=  \frac{f^{m}(s, z- \frac{\theta}{\alpha}, u) - f^{m}(s, z^{'}- \frac{\theta}{\alpha}, u)}{z - z^{'}},  \; \textrm{if}\; z - z^{'} \neq 0,\\
     \\
    \lambda^{m}_{s}(z, z^{'}) &  := \; 0, \quad \textrm{otherwise},\\
    \end{array} \right.
  \]  
the process $\lambda^{m}(Z, Z^{'})$ is in BMO($W$) as soon as both the two processes $Z$ and $Z^{'}$ have this property. Now and for sake of clarity, we denote by $M^{1,m}$  instead of $ Z^{1,m} \cdot W + U^{1, m} \cdot \tilde{N}_{p}$ the martingale part of $Y^{1,m}$. Relying on the relation: $f^{1,m}(s, 0,0) \equiv 0$, we apply the It\^o formula to $Y^{1, m}$ between $t$ and $\tau$ ($\tau$ being an arbitrary stopping time such that: $t \le \tau \le T $) 
 \[ \begin{array}{l}
Y_{t}^{1, m} -  Y_{\tau}^{1, m}  = \\
\quad \quad
\disp{\int_{t}^{\tau}\big(f^{1, m}(s, Z_{s}^{1, m}, U_{s}^{1, m}) - f^{1, m}(s, 0, 0) \big)ds - \big(M_{\tau}^{1,m} - M_{t}^{1,m}\big)}\\
\\
\quad \quad
= \disp{\int_{t}^{\tau}\big(f^{ m}(s, Z_{s}^{1, m} - \frac{\theta_{s}}{\alpha}, U_{s}^{1, m}) - f^{ m}(s, - \frac{\theta_{s}}{\alpha}, 0) \big)ds  - \big(M_{\tau}^{1,m} - M_{t}^{1,m}\big)}\\
\\
    \end{array}  \]
and we then use the following upper bound 
$$ f^{ m}(s, Z_{s}^{1, m} - \frac{\theta_{s}}{\alpha}, U_{s}^{1, m}) - f^{ m}(s, - \frac{\theta_{s}}{\alpha}, 0) \le   Z_{s}^{1, m} \lambda_{s}^{m}(Z_{s}^{1, m}, 0) + \disp{\int_{\mb{R}\setminus\{0\}}U_{s}^{1, m}(x)\gamma_{s}(U_{s}^{1, m}(x), 0)n(dx)}$$
Defining the measure $\mb{Q}^{m}$ by setting: $\frac{d\mb{Q}^{m}}{d\mb{P}}: = \mc{E}_{T}(\lambda^{m}(Z_{s}^{1, m}, 0) \cdot W + \gamma \cdot \tilde{N}_{p}) $, Girsanov's theorem yields that $W^{\lambda^{m}} := W - \lan{\lambda^{m} \cdot W}, W  \ran $ and $ \tilde{N}^{\gamma}(ds, dx) = \tilde{N}_{p}(ds, dx) - \disp{\int_{\mb{R}\setminus\{0\}}\gamma_{s}(U_{s}^{1, m}(x), 0)n(dx)ds}$ are local martingale under $\mb{Q}^{m}$
and $Y^{1, m}$ is the sum of a local martingale and an increasing process. Using a standard localization procedure, there exists a sequence ($\tau^{n, m}$) converging to $T$, as $n$ goes to $\infty$ and such that 
$$ Y_{t}^{1, m} \le \mb{E}_{\mb{Q}^{m}}\left( Y_{\tau^{n,m}}^{1, m}| \mc{F}_{t}\right),$$
and inequality (\ref{eq: estimessent}) follows from the application of the bounded convergence theorem to
($\mb{E}_{\mb{Q}^{m}}\left( \tilde{Y}_{\tau^{n,m}}| \mc{F}_{t}\right)$)$_{n}$ and the almost sure convergence of $ \tilde{Y}_{\tau^{n, m}}$ to $\frac{B}{N},$ resulting from the fact that $(\tau^{n, m})_{n}$ becomes stationnary, $\mb{P}$-a.s.\\ 

\indent To obtain the lower bound, i.e. $Y^{1,m} \ge -\frac{|B|_{\infty}}{N}$, we apply the same procedure to $\bar{Y}^{1, m} = - Y^{1, m}$: in this case, this consists in linearizing the increments of
$$ -f^{1, m}(s, Z_{s}^{1, m}, U_{s}^{1, m}) = -f^{1, m}(s, Z_{s}^{1, m}, U_{s}^{1, m}) -(-f^{1, m}(s, 0,0)).$$
Hence, provided we replace $\lambda^{m}(Z^{1, m}, 0) $ by $\lambda^{m}(0, Z^{1, m})$ and $ \gamma(U^{1, m},0)$ by $\gamma(0, U^{1, m})$, we obtain the same controls as in ($H_{2}$) and rewritting identically the previous proof, it entails: $-Y_{s}^{1, m} \le \frac{|B|_{\infty}}{N}, \; \mb{P}\textrm{-a.s. and for all} \; s,$
which achieves the proof of (\ref{eq: estimessent}).
\begin{flushright}
 $\square$
\end{flushright}
\subsection{A2: Omitted proof of the first stability result}
We prove here the strong convergence of ($Z^{1, m}$) and ($U^{1, m}$) skipped in Section 3.3.1 and which is the essential ingredient in the proof of the stability result in lemma \ref{convuniform}. In all that proof, $C$ stands for an arbitrary constant which may vary from line to line and depends only on the parameters $|B|_{\infty}$ and $\alpha$.
The proof of this result relies on the same methods and computations as in \cite{mkobylanski} but, contrary to the aforementionned paper, we work here in a discontinuous setting, which brings additionnal difficulties.\\
 $(Y^{1, m})$ being increasing, then, for any pair $m,\; p$, such that $p \le m$, $Y^{1, (m,p)} := Y^{1, m} - Y^{1,p}$ is non negative and bounded by $|2\frac{B}{N}|_{L^{\infty}} \leq 2\frac{M_{B}}{N}$ (this results from Appendix A1). Using assertion (i)(b) in Lemma \ref{estim2}, we deduce
$$|U^{1, m,p}_{s}|_{L^{\infty}(n)} \leq 4\frac{M_{B}}{N} , \quad \mb{P}\textrm{-a.s. and for all} \; s.$$
and applying then It\^o's formula to the process $ |Y^{1, (m,p)}|^{2}$, it yields 
 \begin{tabbing}
 \= $  \mathbb{E}\left(|Y_{0}^{1, (m,p)}|^{2}\right)  -\mathbb{E}\left(|Y_{T}^{1, (m,p)}|^{2}\right) = $ \\
\\
 \>
$  \quad   + \; \mathbb{E}\left(\disp{\int_{0}^{T}2Y_{s}^{1, (m, p)}\big(f^{1, m}(s,Z_{s}^{1, m},U_{s}^{1, m})- f^{1, p}(s,Z_{s}^{1, p},U_{s}^{1, p})\big)ds}\right)$ \\
\\ 
\> $\quad   - \; \mathbb{E}\left(\disp{\int_{0}^{T} |Z_{s}^{1, \;(m, p)}|^{2}ds}\right)- \mathbb{E}\left(\disp{\int_{0}^{T}\int_{\mathbb{R}^{*}}|U_{s}^{1, (m, p)}(x)|^{2} n(dx)ds}\right). \quad \quad \quad \quad \quad (*)$ \\
\end{tabbing}
We then need to give an upper bound to the following difference 
\[ \begin{array}{ll} F^{m,\; p}  & = \; f^{1, m}(s, Z_{s}^{1, m}, U_{s}^{1, m}) - f^{1, p}(s, Z_{s}^{1, p}, U_{s}^{1, p})\\
 & = f^{m}(s, Z_{s}^{1, m} -\frac{\theta_{s}}{\alpha}, U_{s}^{1, m}) - f^{ p}(s, Z_{s}^{1, p}-\frac{\theta_{s}}{\alpha}, U_{s}^{1, p}).\; \\
   \end{array} \]
Since both $f^{m}$ and $f^{p}$ satisfy ($H^{1}$), we have
 \[ f^{m}(s, Z_{s}^{1, m}-\frac{\theta_{s}}{\alpha}, U_{s}^{1, m}) \le \frac{\alpha}{2}|Z_{s}^{1, m}-\frac{\theta_{s}}{\alpha}|^{2} + |U_{s}^{1, m}|_{\alpha},\]
and we rely on the classical inequality: $ab \leq \frac{1}{2}(a^{2} + b^{2})$, to obtain 
\[\exists \; \hat{C} \in L^{1}(ds \otimes d\mb{P}), \;  - f^{p}(s, Z_{s}^{1, p}-\frac{\theta_{s}}{\alpha}, U_{s}^{1, p}) \le \hat{C}_{s} + \frac{\alpha}{4}|Z_{s}^{1, p}-\frac{\theta_{s}}{\alpha}|^{2}. \]
with : $\hat{C} = \frac{|\theta|^{2}}{\alpha}$.
Then, we use the convexity of $z \to |z|^{2} $ and $|\cdot|_{\alpha}$ to write, on the one hand,
\begin{tabbing} 
$\frac{\alpha}{2}|Z_{s}^{1, m}-\frac{\theta_{s}}{\alpha}|^{2} $\=$ \leq \; \frac{\alpha}{2}(|\frac{1}{3}(3 Z_{s}^{1, (m, p)} + 3(Z_{s}^{1, p}-\tilde{Z}_{s}) + 3(\tilde{Z}_{s}-\frac{\theta_{s}}{\alpha})|^{2}) $\\
\> $ $\\
\> $\leq  \; \frac{3 \alpha}{2}(|Z_{s}^{1, (m, p)}|^{2} + |Z_{s}^{1, p}- \tilde{Z}_{s}|^{2} + |\tilde{Z}_{s}-\frac{\theta_{s}}{\alpha}|^{2}), $\\
\end{tabbing}
and similarly
$$ \frac{\alpha}{4}|Z_{s}^{1, p}-\frac{\theta_{s}}{\alpha}|^{2} \leq \frac{\alpha}{2}\big( |Z_{s}^{1, p}- \tilde{Z}_{s}|^{2} + |\tilde{Z}_{s}-\frac{\theta_{s}}{\alpha}|^{2}\big) ,$$
and, on the other hand 
\begin{tabbing} 
$|U_{s}^{1, m}|_{\alpha} $ \= $ =\; |(\frac{3 U_{s}^{1, (m, p)}}{3} +\frac{3( U_{s}^{1, p} - \tilde{U}_{s})}{3}  + \frac{3 \tilde{U}_{s}}{3})|_{\alpha}$,\\
\\

 \> $ \leq \; | U_{s}^{1, (m, p)}|_{3\alpha} + |U_{s}^{1, p} - \tilde{U}_{s}|_{3\alpha}+ |\tilde{U}_{s}|_{3\alpha}$,\\ 
 \\
\> $  \leq \; C\big(| U_{s}^{1, (m, p)}|_{L^{2}(n)}^{2} + |U_{s}^{1, p} - \tilde{U}_{s}|_{L^{2}(n)}^{2} + |\tilde{U}_{s}|_{L^{2}(n)}^{2}\big).$\\
\end{tabbing}
To get the first inequality, we use: $|u|_{3\alpha} = \frac{1}{3}|3u|_{\alpha}$,
 and to obtain the constant $C$ appearing in the second inequality, we rely on the relation (\ref{eq: relationeq})
obtained in section 3.2.
Taking into account all these majorations and putting in the left-hand side all terms containing either $Z^{1, (m,\;p)} $ or $U^{1, (m,\;p)}$,
we rewrite It\^o's formula given by (*) as follows 
\begin{tabbing}
 \= $  \mathbb{E}(|Y_{0}^{1, (m,p)}|^{2})) + \; \mathbb{E}\disp{\int_{0}^{T}\big(1 - 4 \alpha Y_{s}^{1, (m, p)} \big) |Z_{s}^{1, (m, p)}|^{2}ds}$\\
\\ 
\>$ \; + \; \mb{E}\disp{\int_{0}^{T}\big(1 -2 C Y_{s}^{1, (m, p)}\big) |U_{s}^{1, (m, p)}|_{L^{2}(n)}^{2}ds} $\\
\\ 
\> $\quad \leq  \; \mathbb{E}\left(\disp{\int_{0}^{T} 2 Y_{s}^{1, (m, p)}\hat{C}_{s} + \; 4\alpha Y_{s}^{1, (m, p)} \big( |Z_{s}^{1, p} -\tilde{Z}_{s}|^{2}+ |\tilde{Z}_{s}-\frac{\theta_{s}}{\alpha}|^{2}\big) ds}\right)$ \\
\\ 
\> $ \quad \quad + \;\mathbb{E}\left(\disp{\int_{0}^{T} 2 C Y_{s}^{1, (m, p)}\big( |U_{s}^{1, p} -\tilde{U}_{s}|_{L^{2}(n)}^{2}+ |\tilde{U}_{s}|_{L^{2}(n)}^{2}\big) ds}\right).$\\
\end{tabbing}
To justify the passage to the limit in each terms of the right-hand side, as $m$ goes to $+ \infty$, $p$ being fixed, we apply Lebesgue's theorem and, for this, we argue \\
\begin{itemize}
 \item $Y_{s}^{1,(m,p)} \to \big(\tilde{Y}_{s} - Y_{s}^{1, p}\big)$, $\mb{P}$-a.s. and for all $s$, as $m$ goes to $+ \infty$ ($p$ fixed),\\
 \item the processes $|Z^{1, p}|^{2}$, $|U^{1, p}(\cdot)|_{L^{2}(n)}^{2}$, $|\tilde{Z} -\frac{\theta}{\alpha}|^{2}$ 
and $|\tilde{U}(\cdot)|_{L^{2}(n)}^{2}$ are in $L^{1}(ds \otimes d\mb{P})$.\\
\end{itemize}
% $(Y_{s}^{m})$ converge vers $\tilde{Y}_{s}$ ainsi que du fait de l'int\'egrabilit\'e dans $L^{1}(ds \%otimes d\mb{P})$ des processus $|Z^{p}|^{2}$, $|U^{p}|_{L^{2}}^{2}$, $|\tilde{Z}|^{2}$ et %$|\tilde{U}|_{L^{2}}^{2}$.
Focusing our attention on the passage to the limit inf, as $m$ goes to $\infty$ ($p$ being always fixed), we use the a priori estimate  
$$\forall \; m \geq p, \quad 0 \leq Y_{s}^{1, (m, p)} \leq 2\frac{M_{B}}{N}, \; \mb{P}\textrm{-a.s. and for all}\; s,$$
and we provide sufficient conditions so that the following terms
\[ \textrm{and} \;\left\{ \begin{array}{l}
  \big(1 - 4 \alpha Y_{s}^{1, (m, p)}). \\
    \\
  \big(1 -2 C Y_{s}^{1, (m, p)}\big).  \\
   \end{array} \right. \]
 are (almost surely) strictly positive. 
This holds as soon as 
$$  (1 - 16\alpha \frac{M_{B}}{N}) \geq \frac{1}{2}\; \textrm{and} \; (1 - 8 C\frac{M_{B}}{N})\geq \frac{1}{2},$$
which provides a constraint condition on $N$ denoted by (\ref{eq: conditiontermF}): under this condition, the two last terms in the left-hand side of It\^o's formula are positive and we obtain 
 \begin{eqnarray*}
\displaystyle{\liminf_{m \to \infty}\mathbb{E}\int_{0}^{T} (1 -4 \alpha Y_{s}^{1, (m, p)})|Z_{s}^{1, (m, p)}|^{2} ds}\\
\quad \quad \geq  \;\disp{\mathbb{E}\left(\int_{0}^{T}(1 - 4\alpha (\tilde{Y}_{s} -Y_{s}^{1, p})) |\tilde{Z}_{s}- Z_{s}^{1, p}|^{2} ds\right)}, \\  \end{eqnarray*} 
and also
 \begin{eqnarray*}
\displaystyle{\liminf_{m \to \infty}\mb{E}{\int_{0}^{T}\big(1 -4 C Y_{s}^{1, (m, p)} \big) |U_{s}^{1, (m, p)}|_{L^{2}(n)}^{2}ds}}  \\
 \geq \; \disp{\mb{E}\left(\int_{0}^{T}\big(1 -4 C(\tilde{Y}_{s} -Y_{s}^{1, p})\big) |\tilde{U}_{s}- U_{s}^{1, p}|_{L^{2}(n)}^{2}ds\right)}.\\
 \end{eqnarray*} 
Rewritting again It\^o's formula 
\begin{tabbing}
 \= $ \mathbb{E}\big(\phi(\tilde{Y}_{0} -Y_{0}^{1, p}\big) +  \;\mathbb{E}\left(\disp{\int_{0}^{T}\big(1 -4\alpha (\tilde{Y}_{s} -Y_{s}^{1, p}) \big) | \tilde{Z}_{s}- Z_{s}^{1, p}|^{2}ds}\right)$\\
\\ 
\>$ \; + \;\mb{E}\left(\disp{\int_{0}^{T}\big(1 -2 C(\tilde{Y}_{s} -Y_{s}^{1, p}) \big) |\tilde{U}_{s} - U_{s}^{1, p}|_{L^{2}(n)}^{2}ds}\right) $\\
\\ 
\> $\quad \leq \; \mathbb{E}\left(\disp{\int_{0}^{T} 2 (\tilde{Y}_{s} -Y_{s}^{1, p})\hat{C}_{s} + 4\alpha (\tilde{Y}_{s} -Y_{s}^{1, p}) \big( |Z_{s}^{1, p} -\tilde{Z}_{s}|^{2}+ |\tilde{Z}_{s}-\frac{\theta_{s}}{\alpha}|^{2}\big) ds}\right)$ \\
\\ 
\> $ \quad \quad + \;\mathbb{E}\left(\disp{\int_{0}^{T} 2(\tilde{Y}_{s} -Y_{s}^{1, p} ) C \big( |U_{s}^{1, p} -\tilde{U}_{s}|_{L^{2}(n)}^{2}+ |\tilde{U}_{s}|_{L^{2}(n)}^{2}\big) ds}\right).$\\
\end{tabbing}
To proceed with a second passage to the limit (as $p$ goes to $\infty$), we transfer into the left-hand side of the previous and last inequality all terms containing either $|Z_{\cdot}^{1, p} -\tilde{Z}_{\cdot}|^{2}$ or $|U_{\cdot}^{1, p} -\tilde{U}_{\cdot}|_{L^{2}}^{2}$, relying again on the condition (\ref{eq: conditiontermF})
to justify
the passage to the limit. For the right-hand side, the use of Lebesgue's theorem is justified arguing that
\begin{itemize}
 \item the processes $\hat{C}$, $|\tilde{Z} - \frac{\theta}{\alpha}|^{2}$ and $|\tilde{U}|^{2}$ are in  $L^{1}(ds \otimes d\mb{P})$,
 \item $Y_{s}^{1, p} \to \tilde{Y}_{s}$, $ \mb{P}$-a.s and for all $s$.
\end{itemize}
Taking the limit sup over $p$ in the left-hand side of It\^o's formula, it leads to
 \begin{eqnarray*}
\disp{\lim_{p \to \infty} \sup  \frac{1}{2}\left(\mb{E}\int_{0}^{T} |\tilde{Z}_{s} - Z_{s}^{1, p}|^{2}ds +  \mb{E}\int_{0}^{T}|\tilde{U}_{s} - U_{s}^{1, p}|_{L^{2}(n)}^{2} ds\right) } \le 0,\\
   \end{eqnarray*} 
the last inequality being an equality, this ends the proof.
\begin{flushright}
$\square$
\end{flushright}

\subsection{A3: Omitted proof in Section 3.3.2 (the second stability result)}
As in the second Appendix, we prove the strong convergence of $(Z^{2, m}) $ and $( U^{2, m})$ (skipped in section 3.3.2): however, in that case, there is an additional difficulty, since the sequence $(f^{2,m})$ is neither increasing nor decreasing. As before
and for any $(m, p)$, we define $Y^{2, (m,p)}$ by: $Y^{2, (m,p)} :=Y^{2, m} - Y^{2, p}$ and similarly $Z^{2, (m,p)} $ and $U^{2, (m,p)} $. We then apply
 It\^o's formula to $ |Y^{2, (m, p)}|^{2}$ between 0 and $T$ and we take the expectation to obtain\\
\begin{tabbing}
 \= $  \mathbb{E}(|Y_{0}^{2, (m, p)}|^{2})  + \;\mathbb{E}\left(\disp{\int_{0}^{T} |Z_{s}^{2, (m, p)}|^{2}ds}\right)
\;   + \; \mathbb{E}\left(\disp{\int_{0}^{T}\int_{\mathbb{R}^{*}}|U_{s}^{2, (m, p)}(x)|^{2} n(dx)ds}\right)$ \\
\end{tabbing}
  \begin{equation}\label{eq: formuleito2}   \leq  \; \mathbb{E}\left(\disp{\int_{0}^{T}2|Y_{s}^{2, (m, p)}||f^{2, m}(s,Z_{s}^{2, m}, U_{s}^{2, m})- f^{2, p}(s,Z_{s}^{2, p},U_{s}^{2, p})|ds}\right). \quad \quad \quad \quad \quad \quad
  \end{equation}

We then give an upper bound of the following quantity
\[\begin{array}{ll}
  F^{m,\; p} & = \;|f^{2, m}(s, Z_{s}^{2, m}, U_{s}^{2, m})- f^{2, p}(s, Z_{s}^{2, p}, U_{s}^{2, p})|,\\
\\
 & \leq \; |f^{m}(s, Z_{s}^{2, m} + Z_{s}^{1, m} - \frac{\theta_{s}}{\alpha}, U_{s}^{2, m} + U_{s}^{1, m})|\\
\\  
  &\quad \; +\; |f^{p}(s,Z_{s}^{2, p} + Z_{s}^{1, p} - \frac{\theta_{s}}{\alpha}, U_{s}^{2, p} + U_{s}^{1, p})| \\
\\
& \quad \; + \;\;  |f^{m}(s, Z_{s}^{1, m} - \frac{\theta_{s}}{\alpha},  U_{s}^{1, m})|  +\; |f^{p}(s, Z_{s}^{1, p} - \frac{\theta_{s}}{\alpha},  U_{s}^{1, p})|.\\
\end{array} \]
Relying again on the assumption ($H_{1}$) satisfied by any $f^{m}$ (with parameters independent of $m$ or of $p$), we claim, using the estimates of lemma \ref{estim2}, that both processes $Z^{1, m}$ and $Z^{1, p}$ (respectively $U^{1, m}$ and $ U^{1, p}$) are bounded independently (of $m$ and $p$) in $L^{2}(W)$ (respectively in $L^{2}(\tilde{N}_{p})$).
Now, to justify the existence of an integrable random variable $G$ (i.e. $G$ in $L^{1}(ds \otimes d\mb{P})$ which dominates
$$|f^{m}(s, Z_{s}^{1, m} - \frac{\theta_{s}}{\alpha},  U_{s}^{1, m})|  +\; |f^{p}(s, Z_{s}^{1, p} - \frac{\theta_{s}}{\alpha},  U_{s}^{1, p})|, $$
 we refer to the following result (already stated in lemma 2.5, page 569 in \cite{mkobylanski}) 
\begin{proposition}
If ($Z^{m}$)$_{m}$ is a sequence of processes on $[0, T]$ such that
$$\exists M > 0, \;  \disp{\sup_{m} \mb{E}\int_{0}^{T}|Z_{s}^{m}|^{2}ds} \leq M,$$
then, there exists a subsequence ($m_{j}$) such that it satisfies
$$ \disp{\sup_{m \in (m_{j})}|Z^{m}|^{2}} \in L^{1}(ds \otimes d\mb{P}).$$
\end{proposition}
Considering appropriate subsequences of ($|Z^{1, m}|^{2}$) and of ($ |U^{1, m}|^{2}$), we can assume w.l.o.g.
$$\disp{\sup_{m}|Z^{1, m}|^{2}} \; \in L^{1}(ds \otimes d\mb{P}) \;\textrm{and} \; \disp{\sup_{m}|U^{1, m}|_{L^{2}(n)}^{2}} \; \in L^{1}(ds \otimes d\mb{P})$$. Besides, since $\frac{|\theta|^{2}}{\alpha} $ is in $L^{1}(ds \otimes d\mb{P})$) ($\theta$ is bounded), we obtain the existence of a random variable $G$ in $L^{1}(ds \otimes d\mb{P})$ such that
$$|f^{m}(s, Z_{s}^{1, m}- \frac{\theta_{s}}{\alpha},  U_{s}^{1, m})|  +\; |f^{p}(s, Z_{s}^{1, p}- \frac{\theta_{s}}{\alpha},  U_{s}^{1, p})| \leq G.$$
We now use the convexity of both $z \to |z|^{2}$ and $|\cdot|_{\alpha}$ to obtain, on the one hand\\
\begin{tabbing} 
$\frac{\alpha}{2}|Z_{s}^{2, m} + Z_{s}^{1, m}- \frac{\theta_{s}}{\alpha}|^{2} $
\=$\leq  \; \frac{3 \alpha}{2}(|Z_{s}^{2, (m, p)}|^{2} + |Z_{s}^{2, p}- \tilde{Z}_{s}^{2}|^{2} + |\tilde{Z}_{2, s}+ Z_{s}^{1, m}- \frac{\theta_{s}}{\alpha}|^{2}), $\\
\end{tabbing}
and, on the other hand,
\begin{tabbing} 
$| U_{s}^{2, m} + \; U_{s}^{1, m}|_{\alpha} $ \=  
 $ \leq \; | U_{s}^{2, (m, p)}|_{3\alpha} + |U_{s}^{2, p} - \tilde{U}_{s}^{2}|_{3\alpha}+ |\tilde{U}^{2}_{ s}+ U_{s}^{1, m} |_{3\alpha}$ \\
\\ 
\> $  \leq \; C\big(| U_{s}^{2, (m, p)}|_{L^{2}}^{2} + |U_{s}^{2, p} - \tilde{U}_{ s}^{2}|_{L^{2}}^{2} + |\tilde{U}_{s}^{2} + U_{s}^{1, m} |_{L^{2}}^{2}\big).$\\
\end{tabbing}
(In the last inequality, the existence of the constant $C$ results directly from the relation (\ref{eq: relationeq}) and using that the two processes $U_{s}^{1,m}$ and $U_{s}^{2, m} $ are in $( L^{2} \cap L^{\infty})(n)$, $\mb{P}$-a.s. and for all $s$).\\
Similarly, we obtain
\[ \textrm{and} \;\left\{ \begin{array}{l}
 \frac{\alpha}{2}|Z_{s}^{2, p} + \;Z_{s}^{1, p} - \frac{\theta_{s}}{\alpha}|^{2} \leq \alpha\big(|Z_{s}^{2, p}- \tilde{Z}^{2}_{ s}|^{2} + |\tilde{Z}^{2}_{s}+ \;Z_{s}^{1, p}- \frac{\theta_{s}}{\alpha}|^{2} \big) \\
\\
|U_{s}^{2, p} + U_{s}^{1, p}|_{\alpha} \leq  C\big(|U_{s}^{2, p} - \tilde{U}^{2}_{s}|_{L^{2}}^{2} + |\tilde{U}^{2}_{ s} + U_{s}^{1, p}|_{L^{2}}^{2}\big),\\
\end{array} \right.\]
%Toutefois, on commet ici un l\'eger abus, car la constante $C_{\alpha, M_{F}}$ dans la majoration pr\'ec\'edente %n'est pas la m\^eme (cette contante d\'ependant de $\alpha$). Il suffit alors de prendre le maximum des deux %constantes.\\ 
 which entails 
 \begin{tabbing}
 $  F^{m,\; p}  $ \= $ \leq G +  \frac{3\alpha}{2}|Z_{s}^{2, (m, p)}|^{2} + \frac{5\alpha}{2}(|Z_{s}^{2, p}- \tilde{Z}^{2}_{s}|^{2}+ |\tilde{Z}_{s}^{2} + Z_{s}^{1, m}- \frac{\theta_{s}}{\alpha}|^{2})$  \\
\\
 \> $ \quad + \;C| U_{s}^{2, (m, p)}|_{L^{2}}^{2} + 2C\big(|U_{s}^{2, p} - \tilde{U}_{ s}^{2}|_{L^{2}}^{2} + |\tilde{U}_{ s}^{2} + U_{s}^{1, m} |_{L^{2}}^{2}\big).$\\
\end{tabbing}
To conclude, we proceed analogously to the proof given in Appendix A1 and we just give below the main steps: writing 
again It\^o's formula given by (\ref{eq: formuleito2}) by putting in the left-hand side all the terms containing either $ |Z_{s}^{2, (m, p)}|^{2} $ or $| U_{s}^{2, (m, p)}|_{L^{2}}^{2} $, it gives 
\begin{tabbing}
\= $  \mathbb{E}(|Y_{0}^{2, (m, p)}|^{2})  + \;\mathbb{E}\left(\disp{\int_{0}^{T} \big(1 - 3 \alpha Y_{s}^{2, (m, p)}\big)|Z_{s}^{2, (m, p)}|^{2}ds}\right)$\\

$\;   + \; \mathbb{E}\left(\disp{\int_{0}^{T}\int_{\mathbb{R}^{*}} \big( 1 - 2C Y_{s}^{2, (m, p)}\big)|U_{s}^{2, (m, p)}|^{2}(x) n(dx)ds}\right)$ \\
\\
$\; \; $\> $ \le \mathbb{E}\left(\disp{\int_{0}^{T} 5\alpha Y_{s}^{2, (m, p)}(|Z_{s}^{2, p}- \tilde{Z}^{2}_{s}|^{2}+ |\tilde{Z}_{s}^{2} + Z_{s}^{1, m}- \frac{\theta_{s}}{\alpha}|^{2})ds}\right)$\\
\\
$\; \; $\> $\; +  \mathbb{E}\left(\disp{\int_{0}^{T} 2CY_{s}^{2, (m, p)} \big(|U_{s}^{2, p} - \tilde{U}_{ s}^{2}|_{L^{2}}^{2} + |\tilde{U}_{s}^{2} + U_{s}^{1, m} |_{L^{2}}^{2}\big)ds}\right).$\\
\end{tabbing}
To achieve the strong convergence of both ($Z^{2, m} $) and ($U^{2, m}$), it remains to justify tweo successive passge to the limit: i.e, a first time when $m$ goes to $+ \infty$, $p$ being fixed, and a second one when $p$ goes to $ +\infty$. As in the first appendix and
%the main difficulty consists in justifying the passage to the limit inf with Fatou's lemma  
to ensure the assumption of positiveness of both these two processes (this for any pair $m, p$)
$$ \big(1 - 8 \alpha Y_{s}^{2, (m, p)}\big)  \quad \textrm{and}\quad   \big(1 - 4C Y_{s}^{2, (m, p)}\big),$$
 we also impose the following constraint condition
$$ (1 - 16\alpha \frac{M_{B}}{N^{(2)}}) \geq \frac{1}{2} \; \textrm{and}\; (1 - 12 C\frac{M_{B}}{N^{(2)}}
)\geq \frac{1}{2},$$
or equivalently
\begin{equation}\label{eq: nvcondittermF} 
\frac{M_{B}}{N^{(2)}}  = \disp{\inf \{ \frac{1}{32 \alpha}, \;\frac{1}{24 C} \} }.
\end{equation}
%(remarquons juste que celle ci est un peu plus restrictive que la contrainte (\ref{eq: conditiontermF}) trouv\'ee %lors de la preuve pour la famille d'EDSR ($f^{1, \;m} $, $\frac{F}{N}$)).
Provided these two conditions hold, the same procedure as for the first stability result in Appendix A2 can be rewritten and leads to
$$ \disp{\lim_{m \to \infty}\sup \mathbb{E}\left(\int_{0}^{T} |Z_{s}^{2, m}- \tilde{Z}^{2}_{s}|^{2}  ds + \int_{0}^{T}|U_{s}^{2, p} - \tilde{U}_{s}^{2}|_{L^{2}}^{2}ds \right)} = 0.$$
%which is the desired property.\\
\indent To conclude, we justify that this proof can be rewritten identically at each step $k$, $k \ge 2$, to obtain the strong convergence of ($Z^{k, m}$) and ($U^{k, m}$). In fact, to show this, we argue that, for any solution ($Y^{k, m}, Z^{k, m} , U^{k, m}$) of the BSDE given by ($f^{k, m}, \frac{B}{N}$), $Y^{k, m}$ satisfies: $|Y^{k, m}|_{\mc{S}^{\infty}} \le \frac{|B|_{\infty}}{N}.$
(this estimate can be justified by the same argumentation as in Appendix A1).
Besides, if we replace $(Z^{1, m})$ and $(U^{1, m})$ respectively by $(\bar{Z}^{k-1, m}) $ and $(\bar{U}^{k-1,m})$ in the previous proof and using that these two aforementionned sequences are uniformly bounded in $L^{2}(W) $ and in $L^{2}(\tilde{N}_{p}) $, the same procedure holds and implies the strong convergence of the sequences ($Z^{k, m} $) and ($U^{k, m}$) provided the condition (\ref{eq: nvcondittermF}) is satisfied. 

\begin{flushright}
$\square$
\end{flushright}

\end{document}